\documentclass[11pt]{amsart}
\usepackage{amssymb}
\usepackage{amsmath}
\usepackage{amsthm}
\usepackage{graphicx}
\usepackage{bm}
\usepackage[usenames,dvipsnames]{xcolor}
\usepackage{subcaption}
\usepackage{tikz}
\usepackage{pgfplots} 
\usetikzlibrary{3d,calc}
\usepackage{mwe}

\usepackage[T1]{fontenc}
\usepackage[english]{babel}
\pgfplotsset{compat=newest}

\usepackage{float}
\usepackage{color}
\usepackage{hyperref}
\hypersetup{
	colorlinks,
	citecolor= Brown,
	linkcolor= link,
	urlcolor= link
}

\usepackage{appendix}
\usepackage{array}
\usepackage{footnote}
\usepackage{booktabs}
\usepackage{hhline}
\makesavenoteenv{tabular}

\usepackage{algpseudocode}
\usepackage{algorithm}
\usepackage{caption}

\usepackage{bbm}

\usepackage{enumitem}

\usepackage{textcomp}

\usepackage{scrextend}
\deffootnote{0em}{1.6em}{\enskip}

\usepackage{wrapfig}

\usepackage{tcolorbox}
\tcbuselibrary{skins}

\usetikzlibrary{positioning}

\usepackage[top=2cm,bottom=2cm,left=2.5cm,right=2.5cm]{geometry}
\definecolor{link}{rgb}{0.18,0.25,0.63}
\definecolor{myred}{rgb}{0.7,0.25,0.2}
\definecolor{mygray}{rgb}{0.8,0.8,0.8}

\numberwithin{equation}{section}

\DeclareMathOperator*{\argmin}{arg\,min}
\DeclareMathOperator*{\argmax}{arg\,max}

\newcommand{\cmmnt}[1]{}

\makeatletter
\g@addto@macro{\endabstract}{\@setabstract}
\newcommand{\authorfootnotes}{\renewcommand\thefootnote{\@fnsymbol\c@footnote}}%
\makeatother

\definecolor{myred}{rgb}{0.78,0.20,0.00}

\newtheorem{definition}{Definition}[section]

\newtheorem{remark}[definition]{Remark}

\newcommand*\prox{\mathop{}\!\mathrm{prox}}

\newcommand*\Diff[1]{\mathop{}\!\mathrm{D}}

\begin{document}

 \begin{center}
 \Large
   \textbf{\textsc{Data-driven methods for quantitative imaging}} \par \bigskip \bigskip
   \normalsize
   \textsc{Guozhi  Dong}\textsuperscript{$\,1$},
   \textsc{Moritz Flaschel}\textsuperscript{$\,2$},
  \textsc{Michael Hinterm\"uller}\textsuperscript{$\,2$}, \\[0.3em]
  \textsc{Kostas Papafitsoros}\textsuperscript{$\,3$},
   \textsc{Clemens Sirotenko}\textsuperscript{$\,2$},
    \textsc{Karsten Tabelow}\textsuperscript{$\,2$}
\let\thefootnote\relax\footnote{
\textsuperscript{$1$}School of Mathematics and Statistics, HNP-LAMA, Central South University, Lushan South Road 932, 410083 Changsha, China
}
\let\thefootnote\relax\footnote{
\textsuperscript{$2$}Weierstrass Institute for Applied Analysis and Stochastics (WIAS), Mohrenstrasse 39, 10117 Berlin, Germany}

\let\thefootnote\relax\footnote{
\textsuperscript{$3$}School of Mathematical Sciences, Queen Mary University of London, Mile End Road, E1 4NS, UK}

\let\thefootnote\relax\footnote{
\vspace{-12pt}
\begin{tabbing}
\hspace{3.2pt}Emails: \= \href{mailto:guozhi.dong@csu.edu.cn}{\nolinkurl{guozhi.dong@csu.edu.cn}},
\href{mailto:flaschel@wias-berlin.de}{\nolinkurl{flaschel@wias-berlin.de}},
\href{mailto:Hintermueller@wias-berlin.de}{\nolinkurl{hintermueller@wias-berlin.de}},\\
 \>\href{mailto: k.papafitsoros@qmul.ac.uk}{\nolinkurl{k.papafitsoros@qmul.ac.uk}},
 \href{mailto:sirotenko@wias-berlin.de}{\nolinkurl{sirotenko@wias-berlin.de}},
 \href{mailto:karsten.tabelow@wias-berlin.de}{\nolinkurl{karsten.tabelow@wias-berlin.de}}
  \end{tabbing}
}
\end{center}
\vspace{-0.8cm}

\begin{abstract}In the field of quantitative imaging, the image information at a pixel or voxel in an underlying domain entails crucial information about the imaged matter. This is particularly important in medical imaging applications, such as quantitative Magnetic Resonance Imaging (qMRI), where quantitative maps of biophysical parameters can characterize the imaged tissue and thus lead to more accurate diagnoses. Such quantitative values can also be useful in subsequent, automatized classification tasks in order to discriminate normal from abnormal tissue, for instance. The accurate reconstruction of these quantitative maps is typically achieved by solving two coupled inverse problems which involve a (forward) measurement operator, typically ill-posed, and a physical process that links the wanted quantitative parameters to the reconstructed qualitative image, given some underlying measurement data. In this review, by considering qMRI as a prototypical application, we provide a mathematically-oriented overview on how data-driven approaches can be employed in these inverse problems eventually improving  the reconstruction of the associated quantitative maps.
\end{abstract}

\section{Introduction}\label{sec:introduction}

The combined and sophisticated use of large imaging data and deep neural networks towards  improving image processing and reconstruction algorithms is nowadays ubiquitous. Essentially most of the current state-of-the-art imaging tasks incorporate some (machine) learning component that provides a priori information or structure to the final image reconstruction. This a priori information is  necessary as imaging problems are often ill-posed and available data are typically degraded and incomplete. A typical example is Magnetic Resonance Image (MRI) reconstruction in which the data are formed from noisy subsampled Fourier coefficients of the contrasting tissue image \cite{Wright_1997}.  In the classical setting of variational, i.e. minimization based image reconstruction, prior information is included via a regularization functional. Often the latter represents a model for imposing sparsity of certain image related features, while suppressing adverse noise. An example for such a functional is the total variation (TV) regularization  \cite{ROF92, chambolle1997image, Lustig_2008} which promotes sparsity of discontinuities (often referred to as ``edges'') of reconstructed image intensities while featuring piecewise constant reconstructions. The latter is often an unwanted side-effect of the TV-regularization leading to the development of generalized versions of TV; see, e.g., \cite{TGV, Knoll_2011}. Rather than defining regularizers ``manually'' on a case-by-case basis, in modern data-driven regimes the regularization functional is merely \emph{learned} automatically from reference data. Inspired by an underlying class of image features (contained also in the training data), this typically allows to impose a more sophisticated structure to the final image . For such a technique, several approaches are conceivable. Here we mention only two. Indeed, for once one may ``train'' the regularizer by  using only the given datum to adapt the reconstruction to its paramount features  in the spirit of unsupervised learning  or one uses large datasets of degraded ground truth data (paired or unpaired) in order to learn the desired image distribution and impose this as prior information \cite{arridge_solving_2019, hintermuller_generating_2019,  kamilov_plug-and-play_2023, monga_algorithm_2021, Ruthotto_2021}. These techniques have been intensively developed concerning \emph{qualitative imaging}.  Regarding the latter, the reconstruction values at pixels do not themselves entail any crucial information (e.g. any  diagnostic value), rather the quality of a reconstruction is judged according to its overall image structure and visual quality. This coins the term \emph{structure-informed learning in imaging} for this category of reconstruction methods.

Prior information to image reconstruction can also be imposed by taking advantage of the physical processes governing the data acquisition. Mathematically, the acquisition physics can often be expressed in terms of (partial (P) or ordinary (O)) differential equations (DE) which then govern the reconstruction process. Considering such DE constraints is particularly relevant in \emph{quantitative imaging} where one is interested in the numerical values of biophysical parameters (in physical units) associated with the imaged object as these values provide objective characteristics of the imaged matter or tissue. Among others, quantitative MRI (qMRI) is a prominent example in this imaging category. Indeed, in qMRI the relevant biophysical parameters are the (longitudinal and transversal) relaxation times ($T_1$ and $T_2$) and the proton density ($\rho$); see Subsection \ref{sec:quantitative_MRI} below for an explanation of these parameters. Specific values of these parameters are related to specific tissue types. Modelling the time evolution of tissue magnetization, given the aforementioned biophysical parameters, leads to the Bloch equations, a system of ODEs associated with every voxel of an imaged domain \cite{Bloch_1946, Wright_1997}. Often, however, such a differential equation model is an idealized representation or a simplification of a far more complex underlying physical process. Sometimes the latter may be even unknown in detail or too complicated to be modelled precisely (from first principles). Under such circumstances, experimental data can provide  insights into the underlying physics and may help to (approximately) learn a mathematical model of the involved physical processes. For the latter, nowadays one benefits from the remarkable versatility and approximation properties of deep neural networks (DNNs) \cite{goodfellow2016deep}. Once models of the underlying physical laws have been learned, then these data-driven models are embedded into the image reconstruction process providing further information \cite{DonHinPap22}. For this overall category of methods, we use the term \emph{learning-informed physics in imaging} throughout this review.

Notably as far as imaging applications are concerned, learning-informed physics approaches are significantly less studied than structure-informed learning techniques. In general, the aim of this review paper is to provide a unified presentation of both categories, with a special focus on recent developments of the former one. For the sake of illustration we use (q)MRI as a prototype application.

The remainder of this paper is organized as follows: In Section \ref{sec:model_based}, we describe model-based methods with an emphasis on qualitative and quantitative MRI reconstruction. We recall the classical variational regularization methods for inverse problems for the qualitative MRI problem and physics-integrated methods (e.g. Magnetic Resonance Fingerprinting and related extensions), but still model-based, for the quantitative setting. We proceed with data-driven approaches in Section \ref{sec:data_driven_methods} for both qualitative and quantitative MRI, starting from dictionary approaches to methods incorporating deep neural networks. We also highlight certain statistics-based post-processing techniques for qualitative MR images towards a subsequent improved reconstruction of the quantitative maps. In the same section we also collect recent developments in learning-inform physics in quantitative imaging, with a focus on learning the Bloch solution map for qMRI. The review ends by a section on FAIR data management in mathematical image processing in Section 4.

\section{Model-based methods}\label{sec:model_based}
In order to set the stage for our main focus on quantitative reconstruction, we start this review by discussing classical model-based methods for qualitative imaging.  In particular, we collect here variational methods with a special focus on popular regularization strategies and solution algorithms. In this way, we believe that the similarities and (yet) differences between qualitative and quantitative approaches become apparent. Regarding quantitative methods, in Subsection \ref{sec:quantitative_MRI} we shall see that typical quantitative imaging problems admit a  variational formulation with an objective inspired by a qualitative approach but subject to certain physics-based constraints which involve (as unknowns) the quantitative parameters of interest. Let us mention that qualitative methods are by now rather classical and there already exists a plethora of reviews, see, e.g.,\ \cite{benning_burger_2018, Bredies_2020, Burger2013, chambolle_pock_2016, Cha_Cas_Cre_No_Po_2010, hintermuller_generating_2019, scherzer2009variational}. Hence, we can be brief here.

When considering the quantities of interest (such as image intensities or physical parameters) as functions over an image domain, rather than vector, matrix or tensor quantities with finitely many entries, respectively, relating to, e.g, a pixel grid which represents the (discrete and finite dimensional) image domain, many of the problems presented here can be studied in an infinite dimensional function space setting. Such a framework raises numerous  interesting and challenging mathematical, algorithmic, and numerical questions. However, for the sake of exposition and unification of presentation we formulate problems in a discrete, i.e., finite dimensional, setting (e.g., on a grid of finitely many pixels/voxels) and make specific remarks regarding their function space analogues.

\subsection{Variational methods for imaging and qualitative MRI}\label{sec:variational_methods}

The classical problem in imaging is to compute an approximation $u$ of an underlying ground truth image $u_{true}$, where $u, u_{true}\in X$ with $X$ denoting a set of functions mapping from some Euclidean space to $\mathbb{R}^{m}$.
When $m=1$, often $u(x)$ models the image intensity at a point $x$. In the discrete setting for $u$ this corresponds to the intensity at a specific pixel (in a 2D image domain). In order to make the distinction from quantitative imaging, we already stress here that in the qualitative regime,  the precise values of $u$, whose range is often dictated by some a priorly chosen scale such as, e.g., $u(x)\in [0,1]$, are not of paramount importance. Rather one is merely interested in salient features of the image, like structure, geometry, contrast, or in any other  information that can be used, for instance in diagnostic procedures upon direct visual inspection, e.g., for detecting tumours.  Typically, $u_{true}$ is not directly accessible, but rather through available (measured) data $y$. A common model for $y$ is
\begin{equation}\label{basic_model}
y=Au_{true}+\eta.
\end{equation}
Here $A:X\to Y$ is the forward operator which models some degradation (e.g. by convolution, incompleteness, etc.) or it reflects an image transformation from the image space $X$ into a data space $Y$. For instance in MRI, it denotes a subsampling $P$ of the Fourier coefficients of the original image $u_{true}$ \cite{Lustig_2008}, i.e., $A:=P\mathcal{F}$ with $\mathcal{F}$ the Fourier transform. In contrast, in computerized tomography (CT) it denotes the Radon transform \cite{Natterer_2001}, which collects (typically a limited number of) integrals along lines reflecting X-ray attenuations by different tissues in the body. The variable $\eta$ in \eqref{basic_model} represents a \emph{noise} component and it is modelled as a highly oscillatory function in $Y$ (with zero mean). Noise is often the result of measurement errors, random degradation processes or transmission losses. Note also that while we consider an additive noise term in \eqref{basic_model} , noise may also occur in the measurement process in a nonlinear way. This may give rise to, e.g., multiplicative noise (see \cite{Goodman_1976}) which then requires an according adaptation of \eqref{basic_model}. Unfortunately, the direct inversion of $A$ and, thus, the exact recovery of $u_{true}$ from $y$ is hampered by the presence of noise, but also by the fact that $A$ is often singular or ill-conditioned. These facts render the reconstruction problem an ill-posed inverse problem. A popular way to overcome the latter is by imposing regularization on the reconstruction process via utilizing {\it prior} information on $u$. Following the variational approach to inverse problems, this leads to solving the minimization problem 
\begin{equation}\label{basic_min}
\min_{u\in X} \mathcal{J}(u):=\mathcal{D}(Au,y) + \mathcal{R} (\alpha ; u).
\end{equation}
Above, $\mathcal{D}$ represents a data discrepancy or fidelity term. It is designed to ensure that the reconstruction, i.e., the minimizer of the overall energy $\mathcal{J}$, is close to the data $y$ in a suitable sense. In fact, its proper choice is often mainly dictated by the statistics of the noise, with $\mathcal{D}(Au,y)=\frac{1}{2} \|Au-y\|_{2}^{2}$ being the appropriate form for Gaussian-type noise. In contrast, when salt-and-pepper or random-valued impulse noise is present, respectively, then the use of the $\ell^{1}$ norm is preferable, i.e., $\mathcal{D}(Au,y)=\|Au-y\|_{1}$, and the Kullback-Leibler divergence \cite{Kullback_1951} is the appropriate choice for Poisson noise. The latter is characteristic for tomographic problems, such as  CT or  Positron Emission Tomography (PET). For the sake of exposition and also motivated by our example modality of MRI where the noise can be considered of Gaussian type, in what follows we will use an $\ell^{2}$ fidelity term. The quantity $\mathcal{R}$ denotes the regularization term or prior with respect to the quantity of interest $u$, and $\alpha>0$ represents a scalar or distributed regularization weight. The main purpose of $\mathcal{R}$ is to counteract ill-conditioning (due to properties of $A$) and to filter noise (due to the presence of $\eta$ in \eqref{basic_model}).

In general, the literature on regularization functionals is rather extensive. It went through a significant development especially after the introduction of the Total Variation ($\mathrm{TV}$) regularization in the early 1990s \cite{ROF92}. Its discrete version corresponds to the $\ell_{1}$-norm of the gradient of the image intensity, balanced by a positive regularization parameter $\alpha>0$ against data fidelity:
\begin{equation}\label{TV_discrete}
\mathrm{TV}(u)=\alpha \|\nabla u \|_{1}.
\end{equation}
\begin{figure}[t]
\centerline{
\includegraphics[width=0.25\textwidth]{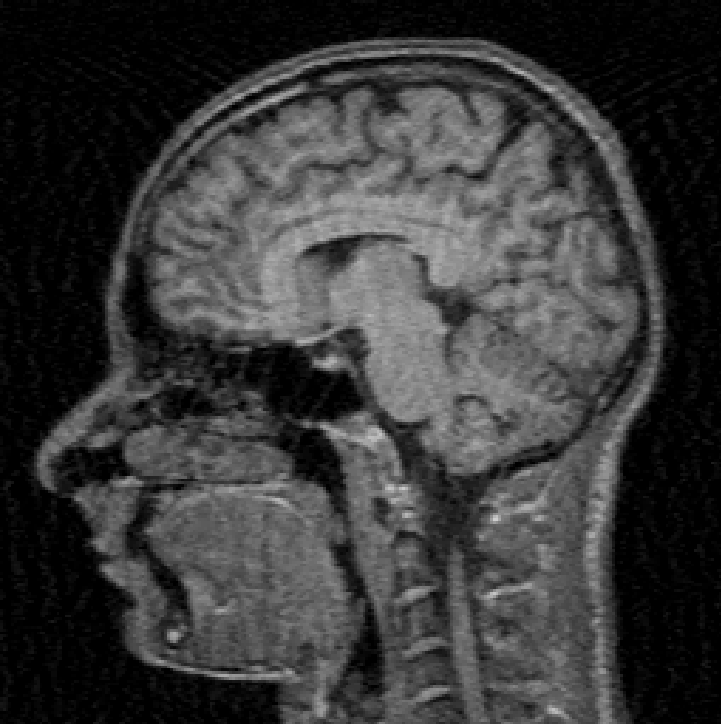} \hspace{0.3cm}
\includegraphics[width=0.25\textwidth]{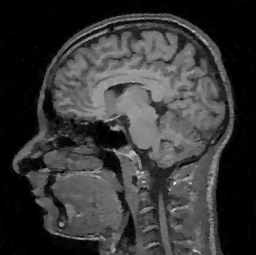} \hspace{0.3cm}
\includegraphics[width=0.25\textwidth]{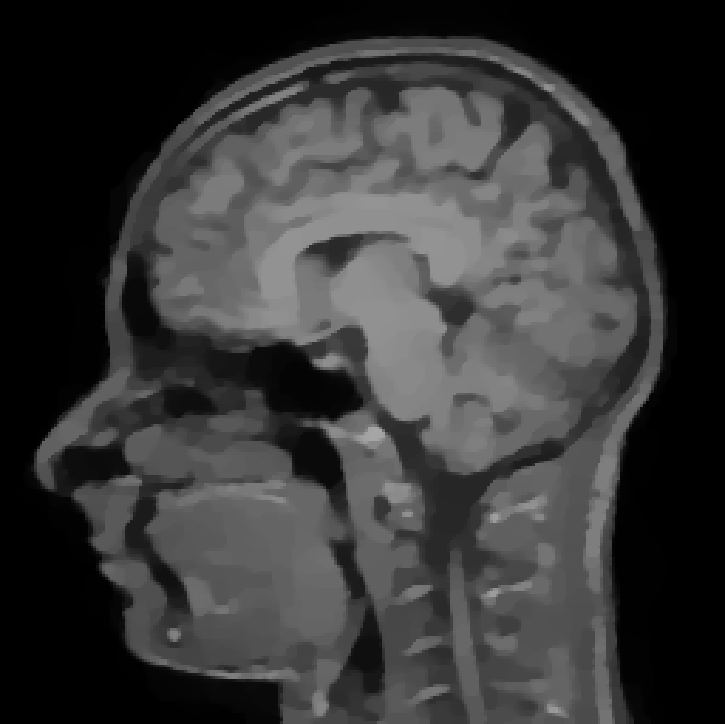}
}
\caption{Qualitative MRI reconstruction example for small (left) and large (right) scalar regularization parameter $\alpha$. Here the data come from brain phantom data (middle) from the Brain Web Simulated Brain Database Data \cite{brainweb,collins1998,kwan1999}.}
\label{mri_over_under_reg}
\end{figure}
TV-regularization based reconstructions typically exhibit a piecewise constant structure, also referred to as \emph{staircasing effect}. The latter is often undesirable. Consequently, higher-order extensions of the total variation have been considered as a remedy. The most prominent choice of such higher order regularizers, is perhaps the Total Generalized Variation (TGV) of order 2 \cite{TGV}, whose discrete version reads 
 \begin{equation}\label{TGV_discrete}
 \mathrm{TGV}(u)=\min_{w} \alpha \|\nabla u- w\|_{1}+\beta \|\mathcal{E}w\|_{1},
 \end{equation}
 where $\alpha, \beta>0$ are regularization parameters and $\mathcal{E}$ denotes the symmetrized gradient $(1/2)(\nabla w+\nabla w^\top)$. It promotes piecewise affine reconstructions, which are often visually more appealing than piecewise constant ones.

\begin{tcolorbox}[
enhanced jigsaw,drop shadow, colback= yellow!75!black, boxsep=0.1cm,
boxrule=1pt, width=1\textwidth, 
interior style={top color=mygray!20!white,
bottom color=mygray!20!white}, 
 opacityback=1,
fonttitle=\bfseries, arc=5pt]
\begin{remark}[Function space versions of classical regularizers]
 In the function space setting, $\mathrm{TV}$ is defined as the total variation of the finite Radon measure $Du$ that represents the distributional derivative of a function $u\in L^{1}(\Omega)$. When indeed the derivative of $u$ admits such a representation, then we say that $u\in \mathrm{BV}(\Omega)$, the space of \emph{functions of bounded variation} \cite{MR0775682,AmbrosioBV}. Similar are the corresponding definitions for TGV.
The study of TV, TGV and other derivative-based regularizers  in their respective function space formulation has also received considerable attention in imaging. This is due to their useful mathematical properties that can provide information about the structure imposed on images, such as the preservation of edges (discontinuities in $u$) and the promotion of piecewise constant/affine structures \cite{MR1852741,caselles2007discontinuity,structuralTV,scherzer2009variational, ring2000structural,papafitsoros2013study}. 
\end{remark}
\end{tcolorbox}

Very large values of regularization parameters typically result in over-regularization leading to over-smoothed reconstructions. While this might benefit the reconstruction of large homogeneous areas in images (where $\nabla u\approx 0$) it can also result in cartoon-like reconstructions with a significant loss of details. On the other hand, very small values of regularization parameters generally preserve details but may also yield reconstructions with noticeable noise and artifacts due to the underlying ill-posedness; compare Figure \ref{mri_over_under_reg} for a representative example in (qualitative) MRI reconstruction. This suggests a local adaptation of the regularization strength which can be achieved by a \emph{spatially varying regularization parameter / weight}. Here the parameter $\alpha$ (and also $\beta$ in TGV) is spatially varying (pixel-dependent) rather than a scalar only. Thus, one has $\alpha = \alpha(x)$, and similarly for $\beta$. Then the corresponding regularization functionals, i.e., \emph{weighted} TV and TGV,  are defined in their discrete version respectively as follows:
\begin{align}
 \mathrm{TV}(u;\alpha)
 &=\|\alpha \nabla u\|_{1},\label{TV_spatial_alpha}\\
  \mathrm{TGV}(u;\alpha,\beta)
 &=\min_{w}\|\alpha (\nabla u-w)\|_{1}+ \|\beta \mathcal{E}w\|_{1}.\label{TGV_spatial_alphabeta}
\end{align}
Clearly, the automatic determination of the spatially varying parameters $\alpha$ and $\beta$ along with $u$ represents a formidable challenge: (i) One has to identify criteria for the quality of a choice of $(\alpha,\beta)$. Such criteria will depend on $u$ and the noise statistics. (ii) A mathematical model for jointly finding $u$ and $(\alpha,\beta)$ has to be developed. The formulation needs to accomodate the identification of several quantities which enter the problem in different ways. (iii) An efficient algorithmic treatment has to be developed to keep the overall computational burden acceptable. Indeed, one has to cope with a non-convex and also non-smooth optimization task.

\begin{figure}[t]
\centerline{
\includegraphics[width=0.25\textwidth]{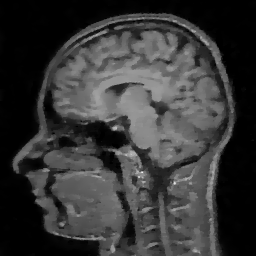} \hspace{0.3cm}
\includegraphics[width=0.25\textwidth, trim={0.6cm 0.6cm 0.5cm 0.8cm},clip]{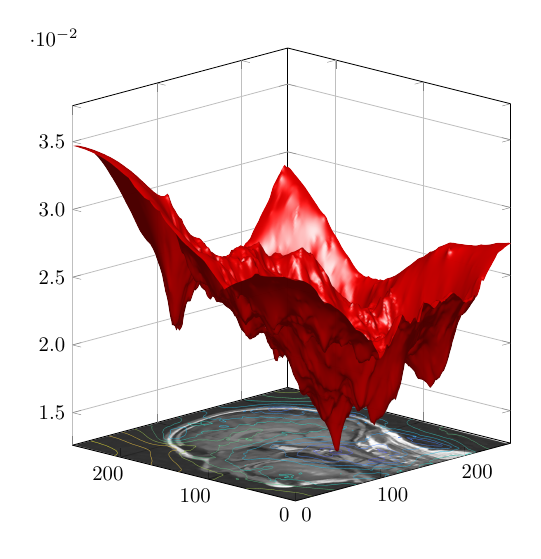} \hspace{0.3cm}
}
\caption{Qualitative MRI reconstruction (left) using TGV with a spatial varying regularization parameter $\alpha$ (right), computed via a bilevel optimization scheme. Observe how the weight attains small values in the bottom right part of the image, where most of the fine scale details are present, leading to their preservation.}
\label{fig:bilevel}
\end{figure}

It has turned out that \emph{primal-dual equilibrium formulations} are particularly suitable for these tasks. In \cite{MR2782120} for TV-regularization and Gaussian noise an augmented Lagrangian type update scheme for a fidelity (\emph{not} regularization) weight was proposed, while simultaneously reconstructing $u$. It utilizes Gumbel statistics, i.e., statistics of the extremes, to develop a local quality measure for the associated fidelity weight. Indeed, the spatially varying weight can be interpreted as a Lagrangian multiplier for a localized version of the image residual. Then the multiplier gets updated as long as the localized image residual does not satisfy a local variance bound; otherwise the fidelity weight remains fixed. The choice of the variance bound is delicate and utilizes the statistics of the extremes. This approach was later extended to $\ell^1$-TV reconstruction \cite{MR2658822}, the reconstruction of color images \cite{MR2780771}, and to TGV \cite{MR3018412}. Another very suitable approach to simultaneously identifying $u$ and $(\alpha,\beta)$ relies on \emph{bilevel optimization} which we briefly describe next, but we also refer the interested reader to the reviews \cite{hintermuller_generating_2019,calatroni2017bilevel,DelosReyes2021}. The bilevel formulation for computing optimal spatially varying regularization parameters reads
\begin{equation}\label{bilevel_general}
\left \{ 
\begin{aligned}
&\min_{u,\mathbf{\alpha}}\;\; \mathcal{L}(u,\mathbf{\alpha})\\
&\text{subject to }\;\;  u\in \argmin_{u\in X}\; \|Au-y\|_{2}^{2} + \mathcal{R} (\mathbf{\alpha} ; u),
\end{aligned}\right.
\end{equation}
where $\mathbf{\alpha}$ denotes all relevant regularization parameters, such as for instance $\alpha$ in \eqref{TV_spatial_alpha} or $(\alpha,\beta)$ in \eqref{TGV_spatial_alphabeta}. Furthermore, $\mathcal{L}$ denotes an upper  level objective which is to be minimized over both the image $u$ and $\mathbf{\alpha}$ under the constraint that $u$ is the solution of the variational problem with regularization parameter $\mathbf{\alpha}$. Borrowing terminology from mathematical game theory, one may consider $\mathbf{\alpha}$ the ``leader'' and $u=u(\mathbf{\alpha})$ the ``follower'', that accepts any choice of $\mathbf{\alpha}$ issued by the leader and reacts optimally by minimizing the objective in the constraints of \eqref{bilevel_general}. Consequently, as in the above primal-dual equilibrium setting $\mathcal{L}$ has to implement criteria for measuring the quality of $u(\mathbf{\alpha})$. For instance, $\mathcal{L}$ may equal the distance of $u(\mathbf{\alpha})$ to the ground truth $u_{true}$, i.e. $\|u(\mathbf{\alpha})-u_{true}\|_{2}^{2}$. In practice, this approach typically involves a training set of true data in the spirit of supervised learning \cite{calatroni2017bilevel, DelosReyes2021} which, however, makes the approach susceptible to overfitting \cite{Reyes_Villacis_SIAM}. This can be remedied by developing ground-truth-free upper level objectives like the one introduced in \cite{HiRa17,MR3712428} for TV, which is inspired by a dual version of the localized residual estimators of \cite{MR2782120}. Later this approach 
was also used for TGV and general convex regularizers \cite{bilevel_TGV, Pagliari_2022}. As in \cite{MR2782120}, there are statistics-based objectives aiming to enforce a localized version of the residuals $Au-y$ that allows the computation of spatially varying regularization parameters without using any ground truths. In  Figure \ref{fig:bilevel} we depict a TGV-based qualitative MRI reconstruction where the regularization parameter $\alpha$ of TGV is spatially varying and computed by the ground-truth-free bilevel approach in \cite{bilevel_TGV, Sirotenko_master}.

\begin{tcolorbox}[
enhanced jigsaw,drop shadow, colback= yellow!75!black, boxsep=0.1cm,
boxrule=1pt, width=1\textwidth, 
interior style={top color=mygray!20!white,
bottom color=mygray!20!white}, 
 opacityback=1,
fonttitle=\bfseries, arc=5pt]
\begin{remark}[Function space versions of bilevel problems in imaging]
Bilevel problems in imaging of the type \eqref{bilevel_general} have also been extensively studied in infinite dimensions; see all of the aforementioned references for bilevel approaches. Here additional challenges arise. For instance, one must ensure that the pairing of the (now function) $\mathbf{\alpha}$ with the gradient of $u$ (now measure $Du$) e.g.\ $\int \alpha(x) d|Du|$ is well defined. The latter is true for instance when $\alpha$ is continuous \cite{HiRa17, Pagliari_2022}, see \cite{davoli2023dyadic} for an investigation concerning even weaker requirements. Additional complications arise in the case of vanishing weights \cite{structuralTV, baldi2001weighted}, which is typically avoided by constraining the corresponding minimization variable  in \eqref{bilevel_general}  to be bounded away from zero. Further challenges exist  when devising numerical solution algorithms for \eqref{bilevel_general}.  In general, these problems fall into the realm of mathematical programs with equilibrium constraints (MPECs); see, e.g., \cite{MR1419501,MR1641213} for rather general accounts in finite dimensions and \cite{kopacka} for a comprehensive discussion of MPEC stationarity for a specific problem class in infinite dimensions. Typically, MPECs suffer from constraint degeneracy. As a consequence, for the derivation of stationarity conditions it requires advanced non-smooth analysis techniques other than the classical Karush-Kuhn-Tucker (KKT) theory; see \cite{Reyes_Villacis_SIAM, kopacka}. More specifically, one is usually confronted with a variety of stationarity conditions; some providing strong conditions avoiding spurious solutions, but being hard to compute, whereas others are somewhat easier to compute, but too weak a filter for (local) optimality. Also, due to the structure of the feasible set, the design of solution algorithms is challenging for this problem class. 
\end{remark}
\end{tcolorbox}

We close this section by mentioning that there is a broad literature on solvers for the TV- resp. TGV-regularized reconstruction problem. However, as in this review we are primarily interested in different reconstruction models, here we do not dwell on discussing solvers in detail. Rather, we only refer to the first-order primal-dual algorithm for non-smooth convex minimization problems in \cite{chambolle2011first}, the second-order-type semismooth Newton method for TV-regularization in \cite{MR2025014,MR2219285}, the duality-based multigrid scheme in \cite{MR2805356}, and the divide-and-conquer algorithm using graph-cuts in \cite{MR2286448}. This selection indicates the versatility in solver developments for variational regularization approaches.

\subsection{Model-based methods for quantitative MRI}\label{sec:quantitative_MRI}

In addition to the single equation \eqref{basic_model} in qualitative imaging, the quantitative imaging problem also involves a (typically nonlinear) operator $e$ that links the image $u_{true}$ to a (bio-)physical quantity of interest $q_{true}$ via
\begin{equation}\label{eq:e}
e(u_{true},q_{true})=0.
\end{equation}
Often,  $e(u,q)=0$ is given by a system of ordinary or partial differential equations modelling an underlying physical process. Consequently, the aim of quantitative imaging is to identify the value $q$ for every pixel (or voxel in the case of MRI) given some measurements $y$ that correspond to the qualitative parameter $u_{true}$. Hence, the resulting general (and idealized) formulation of the reconstruction problem can be stated as follows: Given $y$, estimate $q$ (close to $q_{true}$) and an associated $u$ (close to $u_{true}$) such that
\begin{equation}\label{eq:qIP}
\left\{
 \begin{aligned}
 e(u,q)=0,\\
y=Au+\eta.\\
 \end{aligned}\right.
\end{equation}
In practice, however, due to ill-posedness of this reconstruction task, one can only expect approximate satisfaction of the above equalities.

In what follows, we focus on quantitative MRI as the underlying imaging modality. But let us also  mention that other modalities fit the framework \eqref{eq:qIP}, too. Two specific examples are Photoacoustic Imaging \cite{DingRenVallelian15,ElbauMindrinosScherzer18} and Magnetic Resonance Elastography \cite{MRE2017}, respectively. 
In the case of MRI we recall that the measured data $y$ represent an incomplete (i.e. subsampled via $P$) set of the Fourier coefficients of an emitted signal. The latter encodes information about the total magnetization $m$ of the nuclear spins of hydrogen atoms with (proton) density $\rho$ of the tissue of interest aligned to a strong external magnetic field $B$.
The underlying physical dynamics of $m$ are governed by the Bloch equations, which are a system of ODEs:
\begin{equation}\label{bloch}
 \begin{aligned}
 \frac{\partial m(t)}{\partial t}&=m(t)\times \gamma B(t)-\left (m_{x}(t)/T_{2}, m_{y}(t)/T_{2}, (m_{z}(t)-m_{eq})/T_{1} \right),
 \end{aligned}
\end{equation} 
given some initial (equilibrium) condition.
The relaxation terms in \eqref{bloch} occur due to the excitation of $m$ from equilibrium state with an electromagnetic radio-frequency wave at resonance frequency of the nuclear spins with the gyromagnetic ratio $\gamma$. The excitation usually distorts $m(t)$ from its equilibrium direction by a flip angle $F\!A$ to some initial value $m(0)$, where the relaxation starts. The quantity $m_{eq}$ is the equilibrium magnetization, and $T_1$ and $T_2$ are the relaxation times of the longitudinal and transverse component of $m$ (relative to $B$) with tissue dependent values, respectively. In case of magnetic field inhomogeneities the transverse relaxation time is slightly different and denoted by $T_2^\star$; in order to access $T_2$ a special re-focussing pulse is needed with the scope of so-called spin-echo sequences. The relaxation process releases again energy in the form of electromagnetic waves which can then be recorded by receiver coils as so-called echos after echo times $T\!E$.
Thus, the equation $e(u,q)=0$ of the general setting \eqref{eq:qIP} takes now the form
\begin{equation}\label{eq:Pi}
u=\rho\mathcal{B}(q)
\end{equation}
where $u=\rho m$, $q=(\rho,T_1,T_2)$ and $\mathcal{B}: q\mapsto m$ is the solution map of the Bloch equations \eqref{bloch}. For convenience of notation, by $\Pi: q\mapsto u$ we denote the map defined in \eqref{eq:Pi}, and sometimes in a slight misuse of notation we also refer to this map as the Bloch solution map.

MRI generates its power from the possibility to employ a variety of different sequences, i.e., spin excitation schemes, each focussing on different tissue properties and resulting in a different image contrast. Often, a second spin excitation is performed after repetition time $T\!R$ even if full recovery of the equilibrium state of the magnetization $m(t)$ has not yet been reached. This is typically reflected in the respective signal equations derived as solutions of \eqref{bloch}. Specific choices of the flip angle $F\!A$, the repetition time $T\!R$, and the echo time $T\!E$ enhance the tissue contrast with respect to $T_{1}$, $T_{2}$ (or $T_{2}^\star$), or $\rho$. This leads to $T_{1}$-weighted, $T_{2}$-weighted, proton density-weighted images, respectively. 
\paragraph{Bloch solution maps for specific excitation sequences}
For some specific magnetization-excitation-sequences the Bloch solution map $\Pi$, or a satisfactory approximation thereof, can be stated explicitly. This is the case, for example, for  a \emph{multi-echo Fast Low-Angle Shot (FLASH)} sequence \cite{Haase1986} for  small flip angles $a$. 
In this case, 
 the Bloch map
 is given by the Ernst equation \cite{Ernst1966}. 
By denoting  the echo time and repetition time of this specific sequence by $T\!E$ and $T\!R$, respectively, then we have (see \cite{Helms2008}):
\begin{align}
  \label{eq:FLASHsimple}
  u 
  &= C \cdot \sin a \cdot
               \frac{1 - e^{- R_1 \cdot T\!R}}{1 - \cos a \cdot e^{- R_1 \cdot T\!R}}
               \cdot e^{-R_2^\star \cdot T\!E} ,
\end{align}
where $R_1 = 1/T_1$ and $R_2^\star = 1/T_2^\star$ is the apparent transverse relaxation time.
Here, $C$ is proportional to the equilibrium magnetization $m$ and thus also to the proton density $C = c\cdot \rho$, with $c$ a spatially varying factor related to the detection sensitivity \cite{Weiskopf2013}.
The variation of $T\!E$ within the multi-echo sequence and the flip angle allows for the estimation of  $R_2^\star$, 
$\rho$ and 
 $R_1$\cite{PolzehlTabelow2023}.
In practice, the approach uses two flip angles, one for proton weighted images and one for $T_1$-weighted images, with six to eight echos at different $T\!E$ each. 
Then, \eqref{eq:FLASHsimple} can be re-parametrized with the so-called ESTATICS model~\cite{Weiskopf2014} for estimating the
apparent transverse relaxation time $R_2^\star$ jointly from the two weightings:
\begin{align}
     \label{eq:ESTATICSa}
      u  &= u_{T_{1}} \cdot e^{-R_2^\star \cdot TE}, \qquad \mbox{for the $T_1$-weighted echos,} \\
     \label{eq:ESTATICSb}
      u  &= u_{PD} \cdot e^{-R_2^\star \cdot TE}, \qquad \mbox{for the proton density-weighted echos.} 
\end{align}
It makes use of the fact that the decay constant \(R_2^\star\) of the exponential
signal decay with $T\!E$, cf.~Eq. \eqref{eq:FLASHsimple}, is identical for all weightings.
Thus, we can simply write 
\begin{equation}
   u = \Pi (u_{PD}, u_{T_{1}}, R_2^\star),
\end{equation}
where $\Pi$ can be written as in \eqref{eq:ESTATICSa} and \eqref{eq:ESTATICSb}.

Another popular flip angle sequence pattern is the \emph{Inversion Recovery balanced Steady
State Free Precession (IR-bSSFP)} \cite{scheffler1999pictorial}. Through this choice, the solution of the Bloch equations can be
simulated by solving a discrete linear dynamical system; see for instance \cite{DavPuyVanWia14, DonHinPap19}.

Returning to \eqref{eq:qIP} and recalling the notation for the Bloch solution map, one needs to estimate $q$ where 
\begin{equation}\label{eq:qIP2}
\left\{
 \begin{aligned}
u&=\Pi(q),\\
y&=Au+\eta.\\
 \end{aligned}\right.
\end{equation}
An alternative and perhaps more practical formulation is 
\begin{equation}\label{eq:two_stage}
\left\{
 \begin{aligned}
u&=\Pi(q),\\
u&=\mathrm{Recon}(y),\\
 \end{aligned}\right.
\end{equation}
where ``$\mathrm{Recon}$'' denotes an arbitrary reconstruction scheme for the qualitative MR image $u$ given the $k$-space data $y$. For instance, in view of our previous discussion on classical variational methods, one may consider
\begin{equation}\label{eq:two_stage_example}
\Bigg\{
 \begin{aligned}
u&=\Pi(q),\\
u&=\argmin_{\tilde{u}}\; \frac{1}{2}\|A\tilde{u}-y\|_{2}^{2}+\mathcal{R}(\alpha;\tilde{u}).
 \end{aligned}
\end{equation}
Note that \eqref{eq:two_stage_example} can be interpreted as a coupled equilibrium system (for finding $(u,q)$) when replacing the convex reconstruction problem by its (necessary and sufficient) first-order optimality or Euler-Lagrange condition.
Below we summarize two families of approaches for tackling \eqref{eq:two_stage}, namely the \emph{two-step approaches} and the \emph{integrated-physics} method.

\paragraph{Two-step approaches for qMRI}
As the name suggests, in this family of approaches the qualitative MRI reconstruction problem, i.e. the second problem in \eqref{eq:two_stage}, is solved first to obtain $u$ and then one uses the first equation to obtain the quantitative parameters $q$.
However, due to the presence of noise, artifacts, the fact that $u$ cannot be expected to equal $u_{true }$, as well as due to challenges that arise from inverting $\Pi$, one relaxes the physical law to arrive at a nonlinear regression problem of the type
\begin{equation}\label{eq:two_stage_relax}
\left\{
 \begin{aligned}
&\min_{q}\; \frac{1}{2}\|\Pi (q)-u\|_{2}^{2},\\
&u=\argmin_{\tilde{u}} \;\frac{1}{2}\|A\tilde{u}-y\|_{2}^{2}+\mathcal{R}(\alpha;\tilde{u}).
 \end{aligned}\right.
\end{equation}
We note that different quantitative reconstruction schemes will result from considering \eqref{eq:two_stage_relax} and solving the qualitative problem with different choices of $\mathcal{R}$, including the spatially adapted TV/TGV from Section \ref{sec:variational_methods}. Two particular instances of this technique are highlighted in the next section on data-driven methods.

The original magnetic resonance fingerprinting (MRF) \cite{Ma_etal13} is a popular method also closely related to \eqref{eq:two_stage_relax}. In MRF, a series of $L$ qualitative images $(u^{(1)},\ldots, u^{(L)})$, corresponding to a sufficiently rich excitation sequence with $L$ read-out times, are reconstructed with no regularisation, i.e., with $\mathcal{R}=0$.  This sequence is chosen, such that the Bloch solution map can be explicitly computed (in closed form, e.g. when using (IR-bSSFP)). In an offline phase, a series of Bloch responses $\mathcal{B}(q)$ pertinent to a range of realistic values for $q$ are simulated (fingerprints) and stored in form of a dictionary denoted by $\mathrm{Dic}$. Finally, the time series that results from considering the values in time of the qualitative images at a specific voxel is mapped to its closest fingerprint. Then this voxel is assigned the value of $q$ that corresponds to that fingerprint. We outline this procedure in more detail in Algorithm \ref{alg:MRD} and note that the associated $\mathcal{B}(q)$ refers to the $x,y$-components of the magnetizations.  
\begin{algorithm}[h!]
		\begin{enumerate}
	\item Qualitative reconstructions:\\
	\emph{Reconstruct the series of $L$ qualitative images $u=(u^{(1)},\ldots, u^{(L)})$ by solving
		\[u^{(\ell)}\in \argmin_{u} \|P^{(\ell)} \mathcal{F} u- y^{(\ell)}\|_{2}^{2}\;\;\; \text{ using } \;u^{(\ell)}=\mathcal{F}^{-1}(P^{(\ell)})^\top y^{(\ell)},\quad \ell=1,\ldots,L. \]}
	\item Voxel-wise matching to the pre-computed fingerprints:\\
	\emph{For every  voxel $i=1,\ldots,N$, extract the closest fingerprint to the time-series $u(i)=(u^{\ell}(i))_{\ell=1}^{L}$ according to \[\mathcal{B}(q_{i})= \argmin_{\mathcal{B}(q)\in \mathrm{Dic}} \; \frac{1}{2}\left \| \frac{\mathcal{B}(q)}{\| \mathcal{B}(q)\|_{2}}- u(i)\right\|_{2}^{2}.\]
		Then assign $\left (T_{1}(i), T_{2}(i)\right )$ to the $i$-th voxel and compute its proton density value $\rho_{i}$  as
		\[\rho_{i}=\frac{\|u(i)\|_{2}}{\|\mathcal{B}(q_{i})\|_{2}}.\] }
		\end{enumerate}
	\caption{\text{Magnetic Resonance Fingerprinting (MRF)  for ${q=(\rho, T_{1}, T_{2})}$ \cite{Ma_etal13}}}
\label{alg:MRD}
\end{algorithm}

We note that the qualitative reconstructions in the first step of MRF suffer from severe artifacts due to the typical subsampling and the absence of regularization. As a result, $L$ must be quite large to produce robust results; see \cite{DonHinPap19} for an analysis of the interplay of these two steps. Another drawback of the method is the computational burden of dictionary matching, since the latter is typically very large. 
Improved versions of MRF aiming to tackle some of these challenges have subsequently appeared. For instance, in \cite{DavPuyVanWia14}  a projected Landweber iteration, called BLIP, is proposed. In this approach the qualitative images are projected onto the dictionary in every iteration; see also related extensions in \cite{Golbabaee_2020}. In \cite{MazWeiTalEld18}, low-rank regularization on the qualitative reconstructions $u$ was introduced. In \cite{mcgivney2014svd} the dictionary containing the pre-computed magnetization trajectories was compressed using a low rank approximation in order to reduce the complexity of the second projection step in \eqref{eq:two_stage_relax}. Building on this idea of simplifying the expensive projection step several more sophisticated methods were developed. A brief summary over existing data-driven methods is provided in Section \ref{sec:data_driven_methods} below and for a very comprehensive and systematic review we refer to the recent overview article \cite{shafieizargar2023systematic}.

\paragraph{Integrated-physics approaches for qMRI}
The previous approaches are highly dependent on the fineness of the dictionary. The latter is related to a sufficiently fine discretisation of the quantitative parameter space in order to yield accurate results. This, however, also results in an increased computational effort. In order to avoid the use of a precomputed dictionary, one can revisit \eqref{eq:two_stage} and impose the physical equation $u=\Pi(q)$ as a constraint instead of projecting onto the dictionary. In that case the regularization $\mathcal{R}$ is applied directly to the quantitive map $q$:
\begin{equation}\label{eq:C_Opt}
\left \{
\begin{aligned}
	&\min_{u, q\in  \mathcal{C}_{ad}} \;\frac{1}{2} \| Au-y \|_{2}^2 + \mathcal{R}(\alpha;q),\\
	&\text{subject to }\;\;  u=\Pi(q),
\end{aligned}\right.
\end{equation}
which  leads to the following reduced (integrated physics) formulation
\begin{equation}\label{eq:Re_C_Opt}
	\underset{q\in  \mathcal{C}_{ad}}{\operatorname{min} }\; \frac{1}{2} \| A\circ \Pi(q)-y \|_{2}^2 + \mathcal{R}(\alpha;q).
\end{equation}
Notice that the minimization for $q$ must be performed over a range of \emph{feasible} (i.e.\ realistic) physical values, denoted here by a convex set $\mathcal{C}_{ad}$. The latter encodes a specific range of values depending on the inspected tissue types.
The above optimization without the explicit regularization term $\mathcal{R}$ was proposed \cite{DonHinPap19}. In that work, a projected Levenberg-Marquardt (L-M) iteration was adopted to solve \eqref{eq:Re_C_Opt}, which is a regularized version of a projected Gauss-Newton scheme. We outline this technique in Algorithm \ref{alg:LM_algo}.

\begin{algorithm}
\begin{algorithmic}
\State Given $q_0\in \mathcal{C}_{ad}$ and a sequence $\{\lambda_n\}_{n\in\mathbb{N}}$ of positive real numbers with $\lambda_n\downarrow 0$, iteratively solve the following problems for $n=0,1,2,\ldots$:
\begin{equation}
\begin{aligned}
	\tilde{y}_n& = y- A\circ\Pi(q_n),\\
	h_n&= \argmin_{h} \; \| \Pi^{\prime}(q_n)h -  A^\dagger\tilde{y}_n\|_{2}^2 +\lambda_n \|h\|_{2}^2, \; \\
	q_{n+1}&= P_{\mathcal{C}_{ad}}(q_n + h_{n} ),
\end{aligned}
\end{equation}
\State Terminate the iteration according to a discrepancy principle; see \cite{DonHinPap19} for details.
	\caption{Levenberg-Marquardt iteration for  physics-integrated  qMRI\newline (e.g. Eq. \eqref{eq:C_Opt} with no explicit regularization) \cite{DonHinPap19}}
	\label{alg:LM_algo}
\end{algorithmic}	
\end{algorithm}

Here, $P_{\mathcal{C}_{ad}}$ is the projection operator onto the feasible set $\mathcal{C}_{ad}$. The operator $\Pi^{\prime}$ is the Frech\'et derivative of the parameter-to-solution map of the Bloch equations, and $A^\dagger$ is the generalized inverse of the sub-sampled Fourier transform, i.e. $A^\dagger= \mathcal{F}^{-1}P^\dagger$, where $P^\dagger$ is the zero filling operator, and $\mathcal{F}^{-1}$ is the inverse Fourier transform.
For comparing a variety of algorithms for qMRI, we present a set of examples in the following. Our tests are based on synthetic data from an anatomical brain phantom, publicly available from the Brain Web Simulated Brain Database \cite{brainweb,collins1998,kwan1999}. More details on how to generate this data can be found in \cite{DonHinPap19} or \cite{DavPuyVanWia14}.

\begin{figure}[h!]
	\centering
\includegraphics[width=0.95\textwidth]{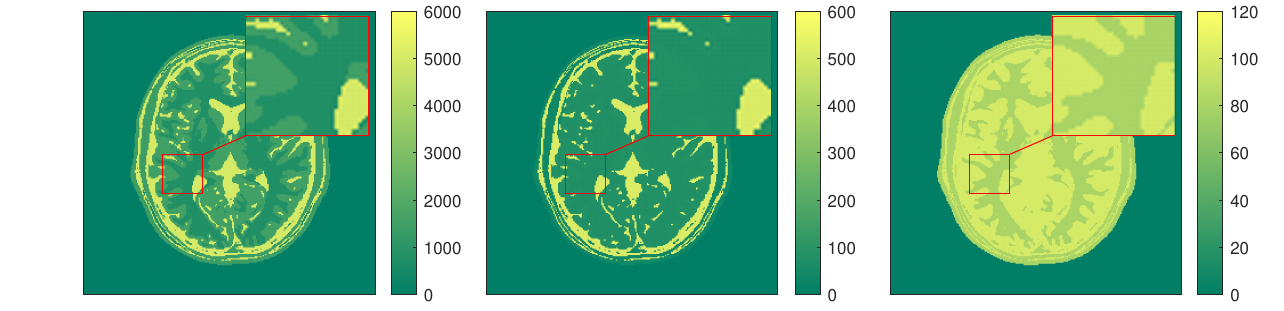}  \\ \vspace{1em}
\includegraphics[width=0.95\textwidth]{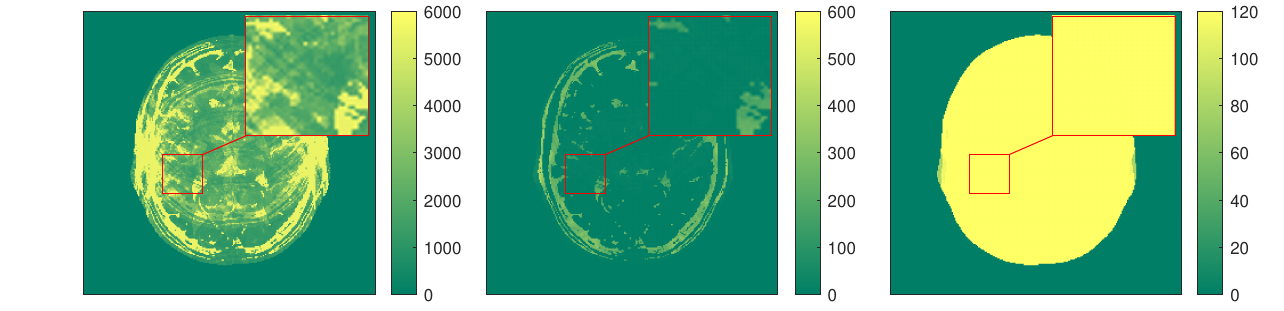} \\ \vspace{1em} 
\includegraphics[width=0.95\textwidth]{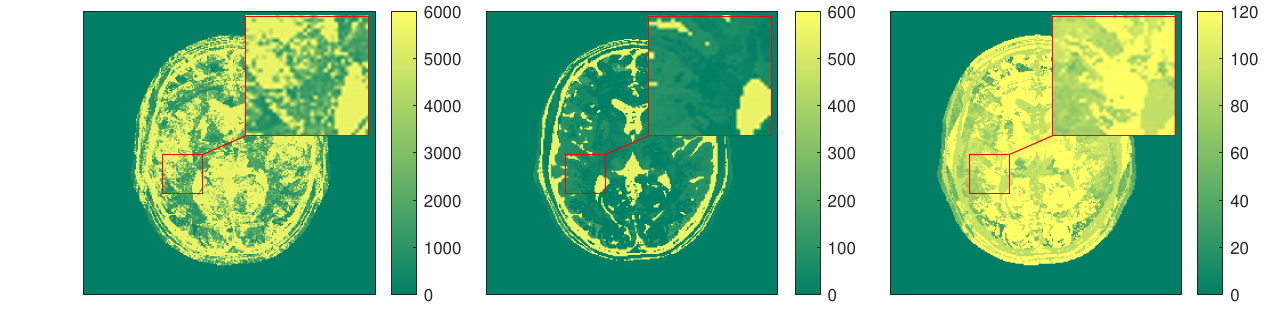} \\  \vspace{1em}
\includegraphics[width=0.95\textwidth]{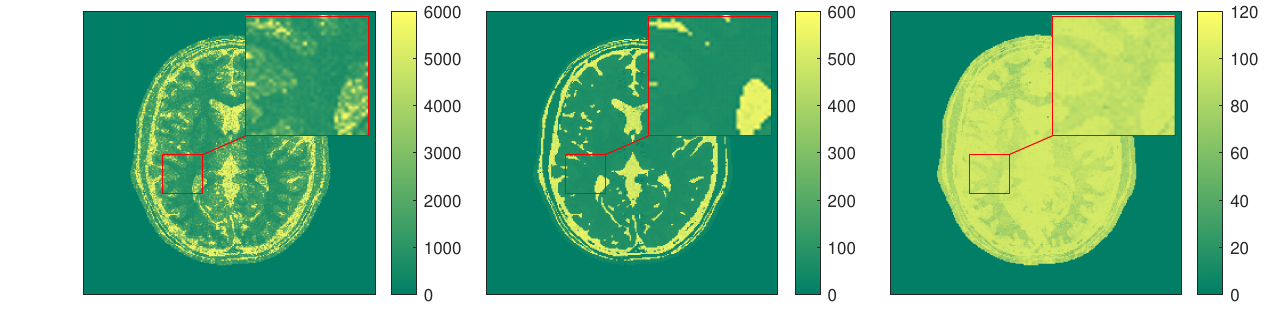} 
\caption{ From left to right (columns): The estimated parameters $T_1$, $T_2$ and the proton density. From top to bottom (rows):  The ground truth, MRF, BLIP, and the integrated model with the L-M algorithm.}
	\label{fig:qMRI_comparisons_parameters}
\end{figure}

\begin{figure}[h]
	\centering
\includegraphics[width=0.95\textwidth]{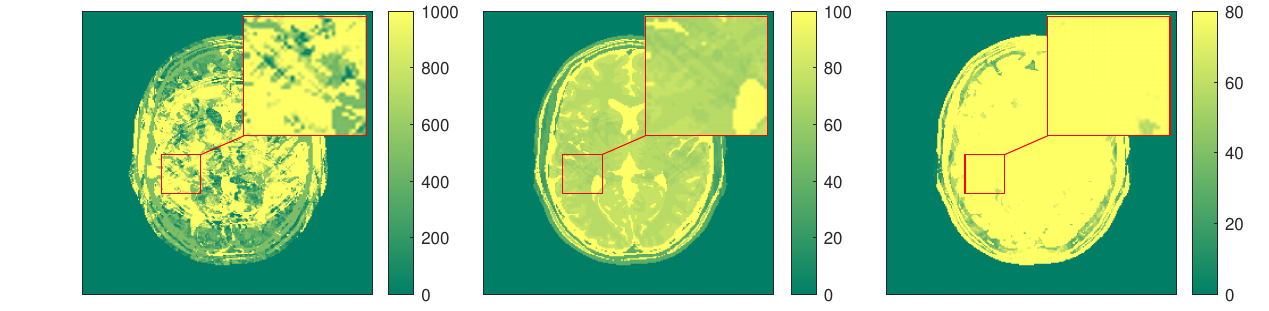} \\ \vspace{1em}
\includegraphics[width=0.95\textwidth]{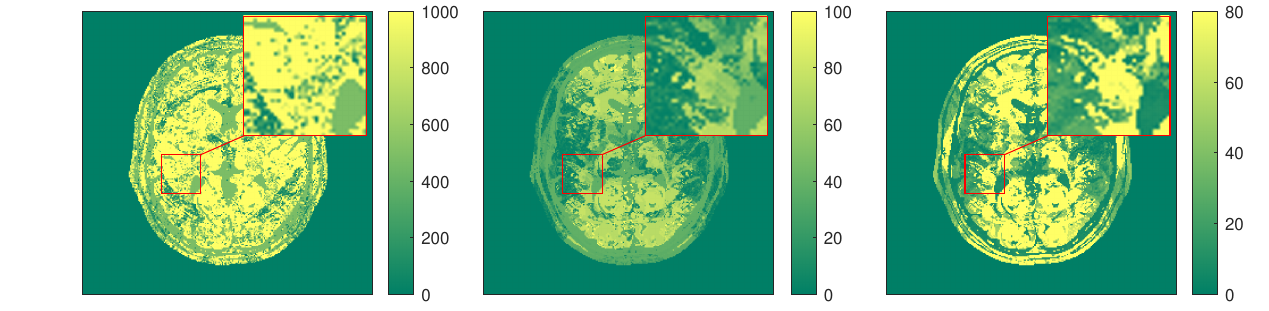} \\ \vspace{1em}
\includegraphics[width=0.95\textwidth]{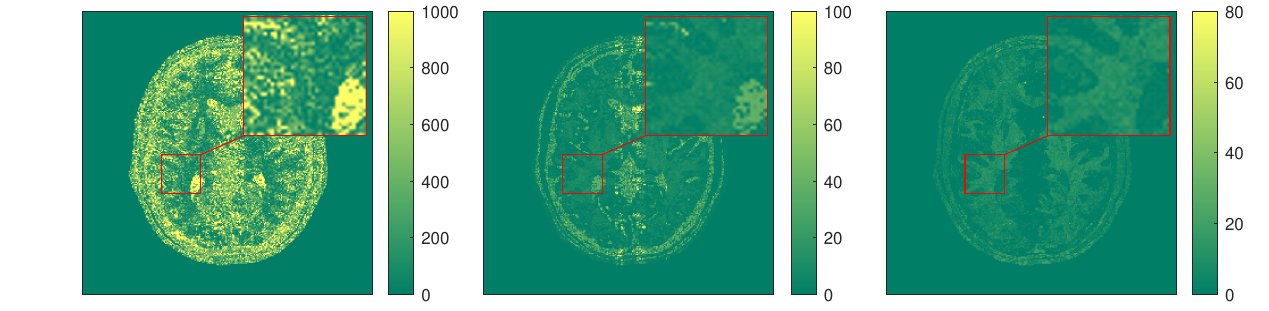} 
	\caption{Error maps of each estimated parameters for  the three algorithms. From left to right: the errors for  $T_1$, $T_2$ and the proton density. From top to bottom:  MRF, BLIP, and the integrated model with the L-M algorithm. Note the different scales in all the three methods reflecting the largest errors in the estimation.}
	\label{fig:qMRI_comparisons_errors}
\end{figure}

The results for each of the above mentioned algorithms can be found in  Figures \ref{fig:qMRI_comparisons_parameters} and \ref{fig:qMRI_comparisons_errors}. We use the MRF reconstruction as initialization for the L-M algorithm, and compare that with the result of the BLIP algorithm when using a relatively refined dictionary. The advantage of integrated-physics approaches is evident by checking the reconstructions in Figure \ref{fig:qMRI_comparisons_parameters} as well as quantitatively by looking at the error maps in Figure \ref{fig:qMRI_comparisons_errors}.
In this example, we have used a time series of the 1/8 Cartesian-subsampled k-space data (Fourier coefficients of magnetizations) of length 40. In the original MRF algorithm (typically requiring a large time series of k-space data), this test setting is far from yielding a reasonable result, while BLIP has improved a little by enforcing the projection to the Bloch manifold. But still one observes deficiencies. The method with integrated physics, however, appears to be efficient and it returns the best results among the three methods.

The results in \cite{DonHinPap19} show that the proposed integrated physics model can work well for settings where an explicit solution formula  for the Bloch equations is available. However, in most cases (excitation sequences) the Bloch equations have no explicit solution formula.  In addition, from a practitioners perspective the Bloch equations appear to be only a simplified mathematical model under some assumptions. As in reality these assumptions are often not realistic, one indeed requires a more complex model.  In this vein, deep learning turns out to be a useful data-driven technique for learning or approximating such physical models. In Subsection \ref{sec:learning_physics} we will therefore discuss extracting physical laws from data.

\section{Data-driven methods}\label{sec:data_driven_methods}
\paragraph{From variational to data-driven regularization}
Traditional regularization techniques as introduced in Section \ref{sec:variational_methods} are based on strong theoretical foundations and often provide good results in many practical applications. However, in many situations important structure in the data may escape ``handcrafted'' regularizers. An example substantiating this claim is depicted in  Figure \ref{fig:wavelet_dec} where the wavelet decomposition of a clean image is shown on the left.  
Natural clean images are assumed to have a sparse wavelet decomposition which is exploited to reconstruct images by solving the problem 
\begin{equation}
\min_{u \in X}\; \frac{1}{2} \| Au-y \|_2^2 + \lambda \|W u \|_1 . \label{eq:wavelet_reg}
\end{equation}
Here $W$ denotes the discrete wavelet transform. Such strategies have been analyzed for instance in
\cite{daubechies2004iterative,grasmair2008sparse,
grasmair2011necessary,candes2006robust,
donoho2006compressed,foucart2013invitation}. Despite good practical results the underlying assumption in \eqref{eq:wavelet_reg} is simply that the wavelet-coefficients of the true image $Wu_{true}$ are sparse in the sense that $\|Wu_{true}\|_0$ is small. However, looking at the first image in  Figure \ref{fig:wavelet_dec}, we clearly observe that image information is ``hidden'' in the coefficients. Indeed, they are not just sparse, but accumulated at edges and highly correlated across different scales. This correlation between wavelet-coefficients is not captured by the simple model in \eqref{eq:wavelet_reg}. To address this issue, one natural strategy is to not only design regularization strategies by taking into account obvious a priori assumptions like sparsity, but to learn parts of the regularization from data. In this vein, during the past decade various methods have emerged; see \cite{arridge_solving_2019} for a comprehensive overview. A more recent article along the same lines is \cite{habring2023neural}, while data-driven reconstruction strategies with a focus on medical imaging are presented in \cite{rueckert2019model,ravishankar2019image}.

One of the first approaches in this direction which directly builds on the ideas of wavelet-regularization is \emph{dictionary learning} where the goal is to find a sparse linear representation of a large amount of given image patches which can then be used as a data-driven regularizer. Dictionary learning usually involves the solution of the non-convex and non-differentiable problem
\begin{equation}
\min_{D \in \mathcal{D}, C \in \mathcal{C}} \|DC - X\|_F^2 + \lambda \|C\|_0,
\end{equation}
where $\mathcal{D}$ is a set of admissible dictionaries, $\mathcal{C}$ is the set of sparse coefficients, $X$ is a collection of clean image patches and again $\lambda>0$ denotes some regularization parameter that controls the sparsity of the coefficients. Here and in the forthcoming part of the article $\| M \|_F$ denotes the classical Frobenius norm of a matrix $M$. In a second step the dictionary can be used as a  regularizer in the sense that the following problem is solved
\begin{equation}
\min_{u,C}\; \frac{\mu}{2} \| Au - y \|_2^2 + \frac{1}{2}\|Ru - DC \|_F^2 + \lambda \|C\|_0,
\end{equation}  
to obtain an estimate for the ground truth.
Here $R$ is an operator that extracts patches from the image, cf. Section \ref{sec:dictionary_learning} for details. It has been recognized later that a decent regularization effect can also be obtained by training the dictionary simultaneously to reconstructing the image $u$. The corresponding method is called \emph{blind compressed sensing} and is reviewed below in Section \ref{sec:dictionary_learning}, where also its applicability to qMRI is demonstrated. Dictionary learning methods in general have been proven to be very successful in medical imaging during the past decade \cite{ravishankar2010mr,li2012group,lingala2013blind}   and are still being used frequently in practice due to their high interpretability, cf. e.g. \cite{pali2021adaptive,ravishankar2019image}. A potential drawback is  the limited expressivity as the final reconstruction method is still of the same  type as \eqref{eq:wavelet_reg} and it is not yet clear whether the approach is suited to capture all complex imaging structures needed to represent natural images. 
\paragraph{Towards interpretable neural networks for linear inverse problems}
For more flexibility and to better model these complex structures, methods based on various types of neural networks have become increasingly popular over the recent years and represent the state of the art in many imaging applications today \cite{lundervold2019overview,wang2020deep,
lin2021artificial}. However, these methods often lead to uninterpretable outcomes and unstable reconstruction algorithms, in the sense that small changes in the setup can lead to significant changes in the image reconstruction \cite{antun2020instabilities,ongie2020deep}.  Nevertheless, the practical results are very impressive, which is why there is a natural demand to combine deep learning strategies with interpretable and robust reconstruction methods for inverse problems. One of the earliest attempts in this direction is called \emph{algorithm-unrolling}. The motivation is to start from an iterative scheme which is known to converge against the solution of \eqref{basic_min} and to interpret this scheme as a neural network. In \cite{gregor2010learning} the ISTA algorithm (see \cite{beck2017first} for details) is under consideration, which is defined subsequently via
\begin{equation}
    u_{k+1} = \prox_{s_k \mathcal{R}(\alpha,\cdot )} (u_k - s_k A^T (A u_k - y)) = \prox_{s_k\mathcal{R}(\alpha,\cdot )} (W_k^1 u_k + W_k^2 y)   \quad u_0 \in X \label{eq:Ista},
\end{equation}
when applied to the problem \eqref{basic_min} with $\mathcal{D}(Au,y) = \frac{1}{2}\| Au - y \|_2^2$. Here we set $W_k^1 = I - s_k A^TA$ and $W_k^2 =  s_k A^T $ and define  $S^k(\alpha,y,W,u_0)$ to be the $k$-th iterate of the algorithm when initialized at $u_0$, with data $y$, regularization parameter $\alpha$, step sizes $s_k>0$ and matrices $W = (W^1_0,\ldots W^1_k,W^2_0,\ldots W^2_k)$. Here $\prox_{s_k \mathcal{R}}:X \to X$ denotes the classical proximal operator from convex analysis, see \cite{rockafellar2009variational}. Note that for $\mathcal{R}(u) = \| u \|_1$ this proximal operator is a piecewise differentiable and nonlinear function which acts component-wise. Hence in this case $S^k(\alpha,y,W,u_0)$ can be understood as a $k$ layer neural network with input $u_0$,  weights $W$ and activation function $\prox_{s_k \mathcal{R}}$. In \cite{gregor2010learning} the $W$-matrices are replaced by general matrices that can be learned from data by solving the following problem
\begin{equation}
    \min_{W} \frac{1}{2M} \sum_{i=1}^M \| S^k(\alpha,y_i,W,u_{0,i}) - u_{true,i} \|_2^2 
\end{equation}
using automatic differentiation such as implemented in nowadays popular software libraries \cite{abadi2016tensorflow, paszke2019pytorch} again using training data
$(u_{true,i},y_i)_{i=1,\ldots,M}$. For an overview on automatic differentiation and how it is applied to learning problems, we refer to \cite{baydin2018automatic}.  It could be shown in experiments that this approach leads to very good outcomes, already when only a few iterations are unrolled. However by replacing all matrices $W^1_k,W^2_k$ by neural networks,  a massive loss of interpretability has to be claimed such that many (open) questions naturally arise: Does the method still converge as the number of unrolled iterates tends to infinity? What is the connection to the optimization problem \eqref{basic_min}, and which regularization properties are inherited from the corresponding variational problem? Moreover due to the slow convergence of proximal gradient type algorithms many iterations need to be unrolled, which might lead to memory issues \cite{hosseini2020dense}. Some research has been conducted in order to overcome these limitations, and also the convergence issue has been addressed, see \cite{chen2018theoretical,liu2019alista} and \cite{chen2022learning} for an overview over existing results. A  natural idea to retain more interpretability compared to pure unrolling is to reduce the number of parameters that are learned in the variational model. For instance, in \cite{ablin2019learning} only the step-sizes are learned. Despite the facts that algorithm unrolling delivers nice practical results in general and also addresses the issues of interpretability to a certain extent there are also  shortcomings, such as, e.g., the large amount of memory consumption. In addition  it is noted in \cite{gilton2021deep} that the method does not converge in all scenarios when $k \to  \infty$, and the reconstruction quality might heavily deteriorate when the number of unrolled parameters is increased during test-time. This is why in \cite{bai2019deep,gilton2021deep} the authors provide a strategy, called (deep) \emph{equilibrium models}, to learn parameters directly in the limit-point of the unrolled method. For this purpose note that any minimizer $u$ of \eqref{basic_min} has to satisfy the fixed point equation
\begin{equation}
    u = \prox_{s \mathcal{R}(\alpha,\cdot )} (u - s A^T (A u - y)) \approx D_\theta (W^1 u + W^2 y),   \quad s>0  \label{eq:Ista_fixedpoint}.
\end{equation}
where again  $W^1 = I - s A^TA$, $W^2 =  s A^T $ and $D_\theta$ is a parametrized approximation of $\prox_{s \mathcal{R}}$, e.g., by a neural network with weights $\theta$. The idea is now to learn $W = (W^1,W^2)$ and $\theta$ from training data by  differentiating through the implicitly given fixed-point in \eqref{eq:Ista_fixedpoint}. The strategy is applied to linear inverse problems in \cite{gilton2021deep} and is further analyzed in \cite{winston2020monotone,bolte2021nonsmooth,
fung2022jfb}. A related method which addresses the convergence problem and parametrizes the prox-operator is the highly flexible Plug-and-Play approach \cite{venkatakrishnan2013plug}. Here also primal-dual algorithms like ADMM  \cite{venkatakrishnan2013plug} or variants of ISTA \cite{gavaskar2020plug} are considered which make use of proximal operators. But in contrast to deep equilibrium models where the parameters $W,\theta$ are learned in an end-to-end fashion, the authors in \cite{venkatakrishnan2013plug} interpret the prox-operator as a denoiser and replace it by a more powerful denoiser $D_\theta :X \to X$, often  parametrized by a neural network with weights $\theta$. In a second step after the denoiser is chosen the method is used iteratively as
\begin{equation}
    u_{k+1} = D_\theta (u_k - s_k A^T(Au_k - y))   \label{eq:pnp},
\end{equation}
where again $s_k>0$ denotes the step size in iteration $k$. After the first attempts of the method using off-the-shelf denoisers like BM3D \cite{dabov2007image} have been proven to be successful, practitioners started to use state-of-the-art methods based on neural networks as denoisers. Due to the possibly complex structure of these denoisers it is a priori unclear whether the plug-and-play algorithm converges or not. This issue has been addressed in \cite{romano2017little,ryu2019plug,gavaskar2020plug,
hurault2021gradient} using quite restrictive classes of denoisers or neural networks, see also the recent overview article \cite{hauptmann2023convergent} where additionally regularization properties of the plug-and-play approach are discussed. 
\paragraph{Data-driven methods for qMRI}
As already pointed out in Section \ref{sec:quantitative_MRI} the existing methods for quantitative MRI can be roughly divided into two different approaches. The first one, explained in  \eqref{eq:two_stage},  reconstructs a high number of MRI-images solving a linear inverse problem to subsequently use these images in a second step to find quantitative parameters that match the reconstructed signal from the first stage. The second class of methods, also called physics-integrated approaches directly computes the physical parameters using the corresponding physical model based on the Bloch equation in  \eqref{bloch}. In principle every data-driven method for linear inverse problems from the introduction of this section can be used to enhance the two-step-procedure. First, the method reconstructs the magnetization trajectories for each pixel in the domain of interest and then the nonlinear projection step is carried out. In \cite{huang2012t2,cao2022optimized,han2022free}  sparsity and low rank regularization techniques are used for the signal reconstruction while in \cite{bilgic2019highly,zibetti2020rapid} deep learning based methods are proposed. A very recent example which also falls into this category and combines total variation based regularization with an unrolling methodology is presented in \cite{kofler2023learning}. Here the starting point is the bilevel problem \eqref{bilevel_general}. In the spirit of the unrolling methodology already discussed above, the lower level problem is replaced by $k \in \mathbb{N}$ steps of an iterative solver (in this case the primal-dual-hybrid-gradient method proposed in \cite{chambolle2011first}) that is known to converge against a solution of the lower level problem. Mathematically speaking, 
\begin{equation}
S^k(\alpha,y,u_0) \to \argmin_{v \in X}\; \frac{1}{2}\|Au - y \|^2_2 + \| \alpha \nabla u \|_1, \quad \text{as } k \to \infty,
\end{equation} 
where $S^k(\alpha,y,u_0)$ denotes the $k$-th iterate of the unrolled algorithm for data $y$, with regularization parameter $\alpha$ and initialization $u_0$. The spatially dependent regularization parameter $\alpha$ is then replaced by a rather small neural network $y \to \alpha_\theta(A^\dagger y)$ with weights $\theta$ and the input being the zero filling solution $A^\dagger y$. The overall learning problem eventually reads
\begin{equation}
\min_{\theta}\; \frac{1}{2M}\sum_{i=1}^M \|S^k(\alpha_\theta(A^\dagger y_i),y_i,A^\dagger y_i) - u_{true,i}\|^2. 
\end{equation}
Here we use again a training set of clean and distorted images $(y_i,u_{true,i})_{i=1,\ldots,M}$ of size $M$. Given new test-data that is generated according to $y = A \circ \Pi (q_{true}) + \eta$ with noise $\eta$ as in \eqref{eq:qIP2}, the method reconstructs first 
\begin{equation}
u^* \in  \argmin_{u \in X}\; \frac{1}{2}\|Au - y \|^2_2 + \| \alpha_{\theta}(A^\dagger y) \nabla u \|_1,
\end{equation}
using $A^\dagger y$ here as the input of the neural network since it already contains significant image structure. Subsequently the nonlinear projection problem in \eqref{eq:two_stage_relax} is solved to obtain an estimate of $q_{true}$. Besides data-driven two-step techniques for qMRI where the data-driven aspect mostly refers to the first linear reconstruction step there are also methods which aim at  modifying the nonlinear projection step in \eqref{eq:two_stage_relax} to infer the tissue parameters. In \cite{mcgivney2014svd} it was empirically shown that the magnetization for different MR sequences actually features a low-rank structure, which was then used to reduce the complexity needed for the tissue parameter inference given a reconstructed magnetization signal. Building on this idea more complex models where used in order to modify or even fully replace the projection step, cf. \cite{asslander2018low,cohen_mr_2018,chen2020high}. Most recent techniques combine data-driven regularization techniques with the first and the second reconstruction step and end up using an end-to-end learning based strategy. We mention here \cite{GOLBABAEE2021101945,song2019hydra} and refer to the comprehensive overview article \cite{shafieizargar2023systematic} for a general account including recent references. 

In the subsequent part of this article we will present three approaches to incorporate data-driven techniques into the reconstruction process for qMRI. The first one, outlined in Section \ref{sec:statistical_methods}, is a two step-approach which uses a data adaptive technique for the linear reconstruction step in \eqref{eq:two_stage_relax}. In particular, it addresses the problem that traditional variational reconstruction methods such as TV reconstruction \cite{ROF92} often lead to a bias in contrast, which might result in a systematic error in the tissue parameters. The second approach, described in Section \ref{sec:dictionary_learning}, uses integrated physics and a regularization strategy based on dictionary learning, which has been proven  successful in particular for unbiased linear estimation of highly under sampled k-space data when no or only a small amount of training data is available. The method which is summarized in Section \ref{sec:learning_physics} falls also into the category of integrated physics but uses training data to modify the solution map of the Bloch equation. In contrast to the neural network based methods described above, the goal here is not only to allow for faster inference of the parameters, but also to account for uncertainties in the physical model.   
\subsection{Data-driven bias-free post-processing methods suitable for qMRI}\label{sec:statistical_methods}

Statistical methods for imaging focus on properties of the error term $\eta$ that is usually considered to be a random variable. 
Its distribution usually depends on the image acquisition method, often, but not always, additive Gaussian noise in equation~\eqref{basic_model} 
is appropriate. 
Specifically, let us assume data $y_i$ at locations $i$ to be distributed as $y_i \sim P_{\theta_i}$ with some density $p(\cdot, \theta_i)$
depending on some local parameter $\theta_i$ (typically from $\mathbb{R}^{p}$) with the probability distribution belonging to some parametric 
(typically exponential) family.

For imaging data usually a structural assumption can be formulated, i.e., there exists a disjoint partitioning of the space of the positions $i$, such that 
$\theta_i$ is constant within each partition. Based on this assumption, an iterative adaptive reconstruction method for the true 
parameter $\theta_i$
called Propagation-Separation approach was defined in \cite{PoSp05}. There, a locally adaptive weighting scheme $W_i^{(k)} =\{ w_{ij}^{(k)}\}$ at a voxel $i$ at iteration
step $k$ is employed for a weighted maximum likelihood estimator
\begin{equation}
    \hat{\theta}_i^{(k)} = \argmax_{\theta}\; \sum_j  w_{ij}^{(k)} \operatorname{log} (p(y_i, \theta_i)),
\end{equation}
where the sum is over all neighboring voxels $j$ with non-vanishing weights $ w_{ij}^{(k)}$. More specifically, the adaptive weights at step $k$ are
defined as 
\begin{equation}
   w_{ij}^{(k)} = K_{loc} (l_{ij}^{(k)}) \cdot K_{st} (s_{ij}^{(k)}),
\end{equation} 
i.e., the product of two kernel functions $K_{loc}$ and $K_{st}$. The first factor in the kernel definition then refers to the Euclidean distance 
\begin{equation}
   l_{ij}^{(k)} = ||i-j||_2 / h^{(k)}
\end{equation} 
in the design space using a bandwidth $h^{(k)}$, the second factor refers to a distance in feature space 
\begin{equation}
   s_{ij}^{(k)} = N_i^{(k-1)} \cdot \mathcal{KL} (\hat{\theta}_i^{(k-1)}, \hat{\theta}_j^{(k-1)}) / \lambda
\end{equation}
where $\mathcal{KL}$ refers to the Kullback-Leibler distance of the two respective probability distributions, $\lambda$ is the adaptation parameter of 
the procedure and $N_i^{(k-1)} = \sum_j w_{ij}^{(k)}$ approximates the variance reduction from the $k-1$-th step.
The resulting procedure, adaptive weights smoothing, with some increasing sequence of bandwidths $h^{(0)}, \dots, h^{(k^\star)}$ performing $k^\star$ steps results in 
noise reduced parameter estimates $\hat{\theta}_i^{(k^\star)}$ that do not suffer from the blurring observed in non-adaptive methods but does 
edge-preserving smoothing. 

\begin{figure}[t]
   \centerline{\includegraphics[width=0.25\textwidth]{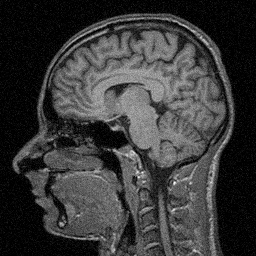} \hspace{0.3cm}\includegraphics[width=0.25\textwidth]{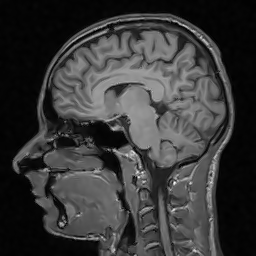} \hspace{0.3cm}\includegraphics[width=0.25\textwidth]{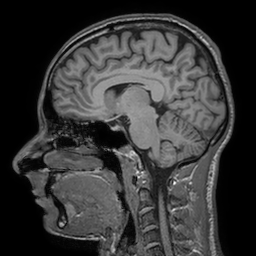} \hspace{0.3cm}}
   \caption{From left to right: MR image with artificial Gaussian noise and fully sampled k-space. Reconstruction with patchwise adaptive smoothing. Original noise-free image. The adaptive smoothing is capable of removing the noise and preserving the details of the image.
   \label{fig:PAWSonMRI}}
\end{figure}

In \cite{PoPaTa20}, a variant of the adaptive weights smoothing employing a statistical penalty $s_{ij}^{(k)}$ as the maximum over a local 
(rectangular) patch of voxels was defined. This patchwise adaptive smoothing combines the iterative procedure outlined above with the idea of 
patchwise comparisons used by non-local means (NLM) \cite{Buades05areview}. The resulting procedure combines the edge preserving property of the original method 
while avoiding the cartoon-like appearance of the reconstructed image that is the result of the structural assumption above. 
In \cite{PoPaTa20}  the method was extensively compared with TV and TGV based reconstruction methods, see the Figure 6 and 7 within this paper. In Figure~\ref{fig:PAWSonMRI},  we provide an
example reconstruction of a noisy MR image (from fully sampled k-space), where $T_{1}$ weighted example data has been taken from the IXI Dataset \cite{IXIdataset}.

\begin{figure}[t]
    \centerline{\includegraphics[width=0.62\textwidth]{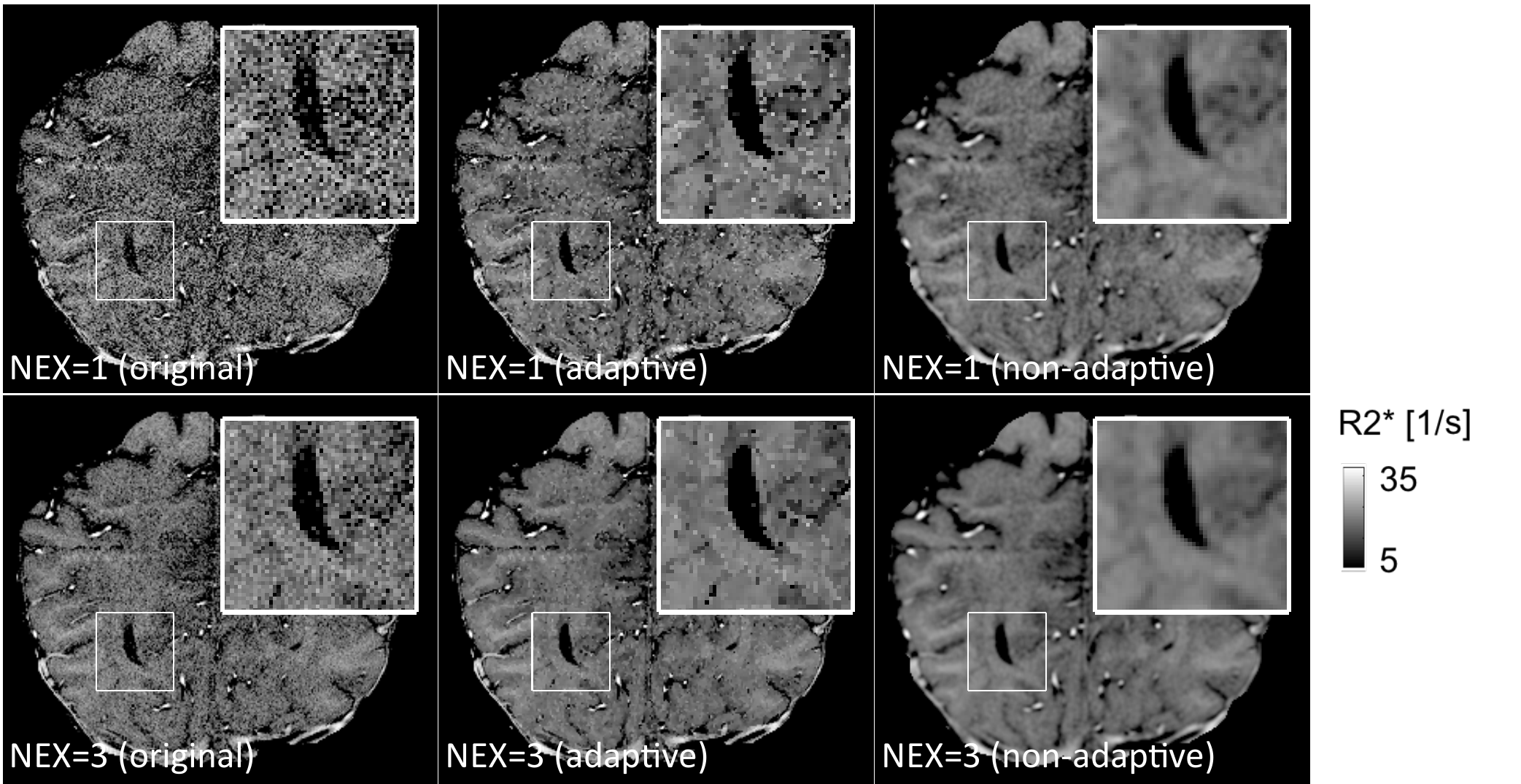} \hspace{0.3cm}}
    \centerline{\includegraphics[width=0.62\textwidth]{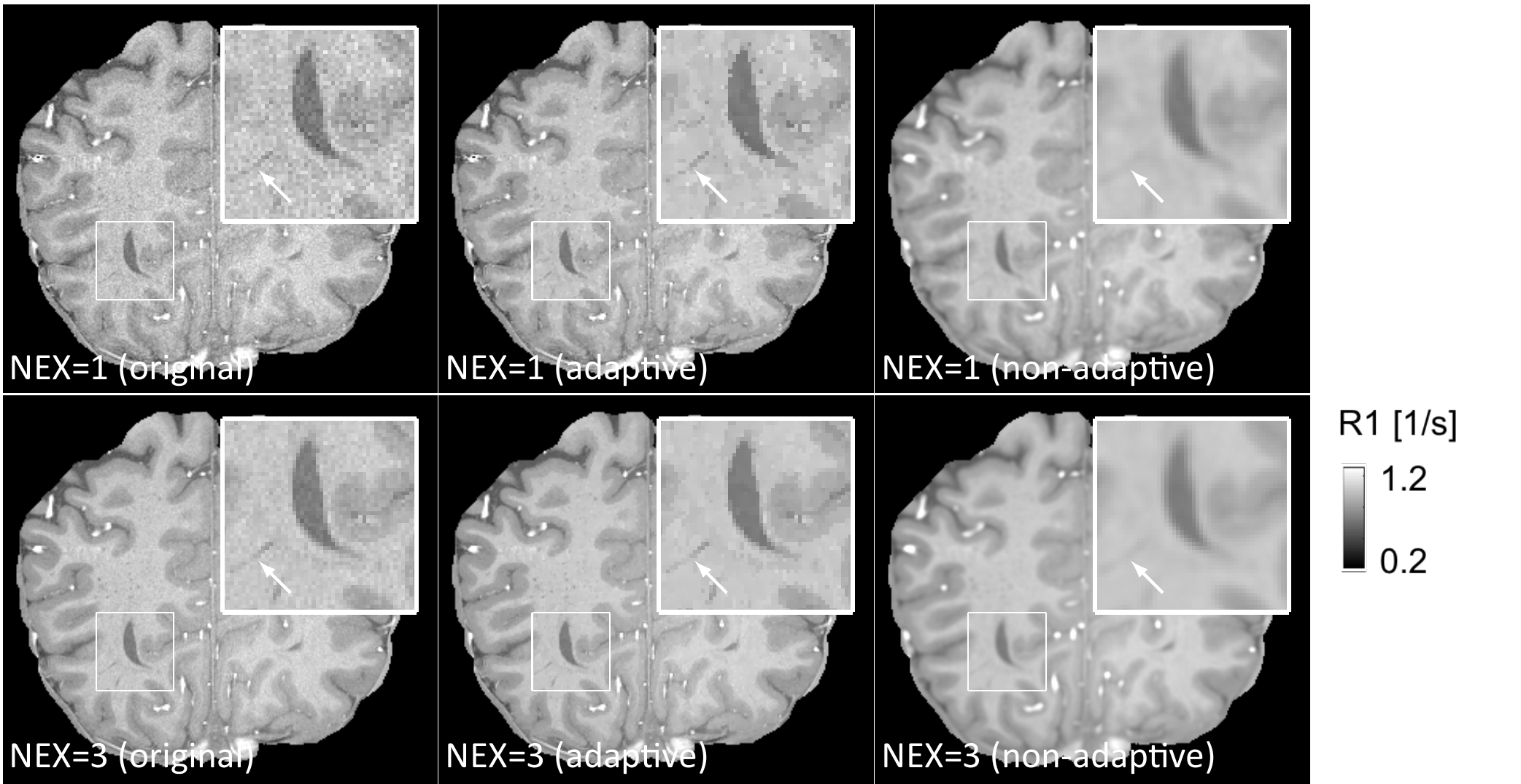} \hspace{0.3cm}}
    \centerline{\includegraphics[width=0.62\textwidth]{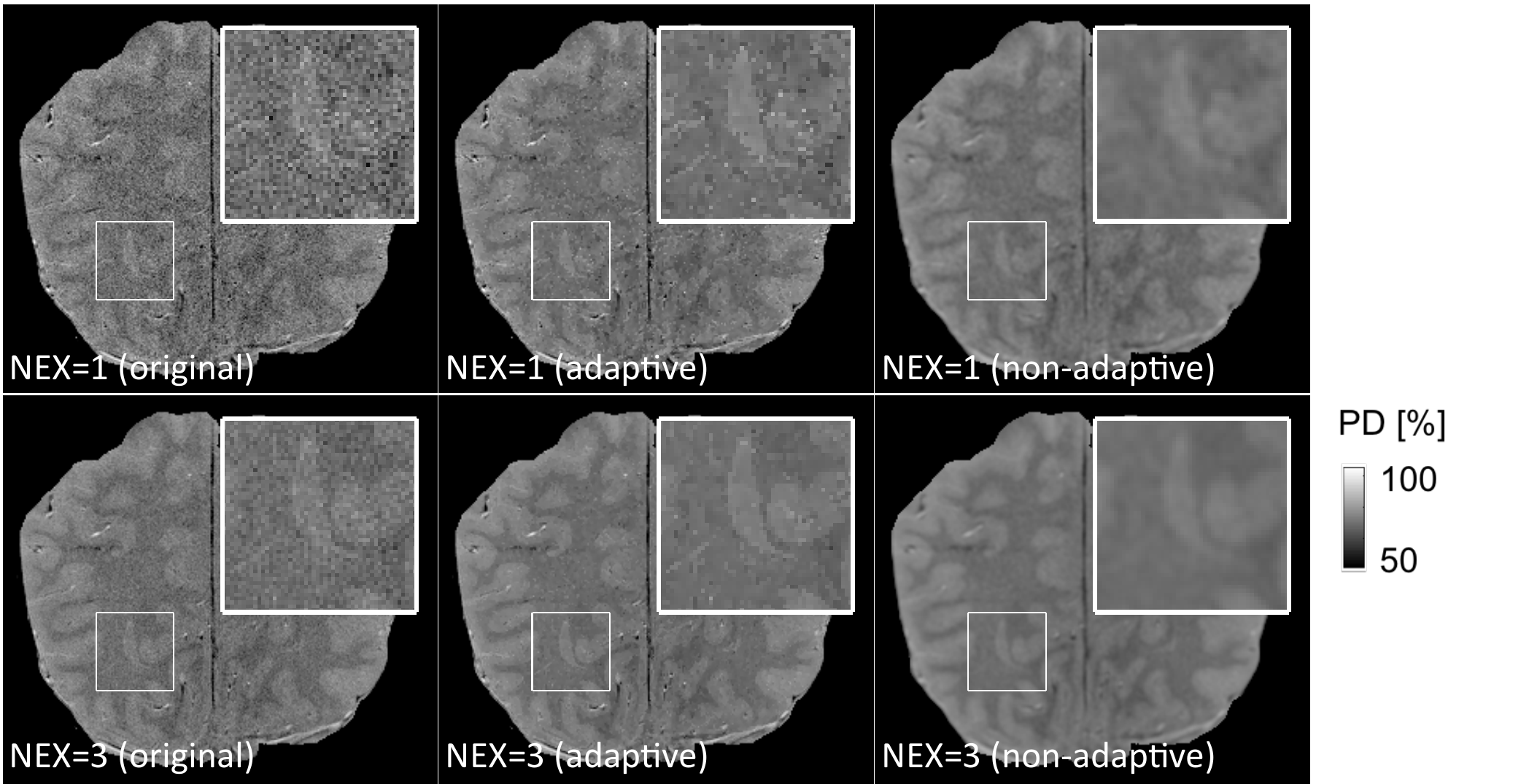} \hspace{0.3cm}}
    \caption{Reconstruction of quantitative MR parameters, i.e., $R_1$, $R_2^\star$, and $\rho$, i.e., proton density PD, after applying adaptive smoothing with comparison to the result obtained from threefold repeated data (NEX=3) and with non-adaptive smoothing, i.e., when a classical Gaussian filtering is applied. While the non-adaptive Gaussian filtering does remove the noise at the cost of blurring, the adaptive method is able to preserve the fine structural details. The comparison between the single data with the threefold repetition shows, that adaptive smoothing is capable of saving image acquisition time by using less data with a comparable resulting map.
    \label{fig:MPMsmoothing}}
\end{figure}

The adaptive weights smoothing procedure can be applied to estimate quantitative MRI parameters using a physical model.  
We will rely on model \eqref{eq:FLASHsimple} and more specifically the ESTATICS re-parametrization  \eqref{eq:ESTATICSa} and \eqref{eq:ESTATICSb} to estimate quantitative parameters like the proton density $\rho$ or relaxation times $T_1$ or $T_2^\star$. First, we solve the optimization problem:
\begin{equation}
  \hat{\theta}_{LS}= \argmin_{\theta}\; \sum_{i=1}^n ({u}_i - \Pi (u_{PD}, u_{T_{1}}, R_2^\star))^2
\end{equation}
to infer on $\theta = (u_{PD}, u_{T_{1}}, R_2^\star)$. The sum is calculated over all $n$ echos of the two sequences from the multi-echo MPM measurement. Comparing \eqref{eq:FLASHsimple} and  \eqref{eq:ESTATICSa}, \eqref{eq:ESTATICSb} it is straightforward to determine the remaining quantitative 
parameters $R_1$, and $A$ (if $T\!R_{T_{1}}$ and $T\!R_{PD}$ are equal) by
\begin{subequations}
    \begin{align}
      \label{eq:R1}
      \hat{R}_1  &= - \ln \left( \frac{\hat{u}_{T_{1}} - \hat{u}_{PD} \cdot \frac{\sin a_{T_{1}}}{\sin a_{PD}}}%
              {\hat{u}_{T1} \cdot \cos a_{T_{1}} - \hat{u}_{PD} \cdot \frac{\sin a_{T_{1}}}{\sin a_{PD}} \cdot \cos a_{PD}} \right) / T\!R,\\
      \label{eq:PD}
     \hat{A}  &= \frac{(1-\cos a_{T_{1}} \cdot e^{-\hat{R}_1 \cdot T\!R}) }{\sin a_{T_{1}} \cdot (1-e^{-\hat{R}_1 \cdot T\!R})} \cdot \hat{u}_{T_{1}}.
    \end{align}
\end{subequations}
The actual measurements within an MPM sequence are corrupted by noise which will be propagated 
to the model parameter maps. The adaptive weights smoothing procedure outlined above can be used 
to reduce the noise in the quantitative maps while preserving the fine structures of the brain tissue that can be seen in the images. 
First we make use of the fact that the ESTATICS\ re-parametrization  \eqref{eq:ESTATICSa} and \eqref{eq:ESTATICSb}
of the signal model~\eqref{eq:FLASHsimple} has a low parameter-induced nonlinearity
and leads to approximate Gaussianity of the estimates
\begin{align*}
  \vec{u}^{(0)} = \left(\hat{u}_{T_{1}}^{(0)}, \hat{u}_{PD}^{(0)}, 
                  \hat{R}_2^{\star {(0)}}\right)^\top.
\end{align*}
A voxelwise estimate $\hat{\Sigma}$ for the covariance of the
parameter estimates can be obtained from the least squares estimation procedure~\cite{PolzehlTabelow2023}.
These three-dimensional parameter maps are then iteratively smoothed by an 
increasing sequence of bandwidths $h_k$ for $k=0, \dots, k^\star$ and the definition of locally adaptive as outlined above.
Specifically, the the statistical penalty $s_{ij}^{(k)}$ defined as
\begin{equation}\label{eq:QMRIStatPen}
  s_{ij}^{(k)} = N_{i}^{(k-1)} \cdot
                 \left(\vec{u}^{(k-1)}_i - \vec{u}^{(k-1)}_j\right) ^\top
                 \hat{\Sigma}_i^{-1}
                 \left(\vec{u}^{(k-1)}_i - \vec{u}^{(k-1)}_j\right)
\end{equation}
based on the estimates $\vec{u}_i^{(k-1)}$ and $\vec{u}_j^{(k-1)}$ and the sum
of weights $N_{i}^{(k-1)} = \sum_j w_{ij}^{(k-1)}$ from the previous step. 
The new estimates in step $k$ are obtained by \begin{equation}
  \vec{u}^{(k)}_i = \sum_j w_{ij}^{(k)} \vec{u}^{(0)} / \sum_j w_{ij}^{(k)}.
\end{equation}
At the last iteration step $k=k^\star$ the final smoothed maps
\begin{align*}
  \vec{u}^{(k^\star)} = \left(\hat{u}_{T_{1}}^{(k^\star)}, \hat{u}_{PD}^{(k^\star)},
                         \hat{R}_2^{\star (k^\star)}\right)^\top,
\end{align*}
can be used within \eqref{eq:R1} and \eqref{eq:PD} to get implicitly smoothed quantitative maps. The application of this procedure 
to real MPM data from \cite{Tabelow2017a} is shown in Figure~\ref{fig:MPMsmoothing}.

\subsection{Dictionary learning approaches for qMRI}\label{sec:dictionary_learning}

\begin{figure}[b]\vspace{2pt}
\begin{center}
\includegraphics[width=0.27\textwidth]{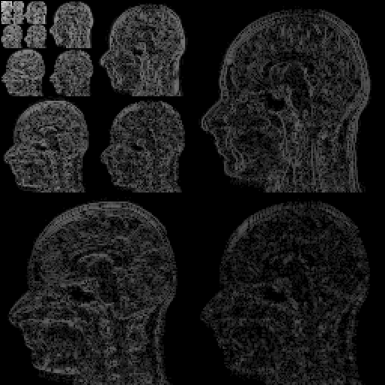}
\hspace{3pt}
\includegraphics[width=0.27\textwidth]{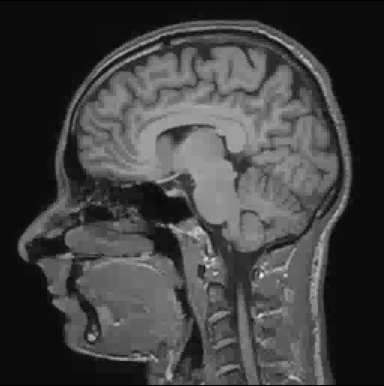}
\hspace{3pt}
\includegraphics[width=0.27\textwidth]{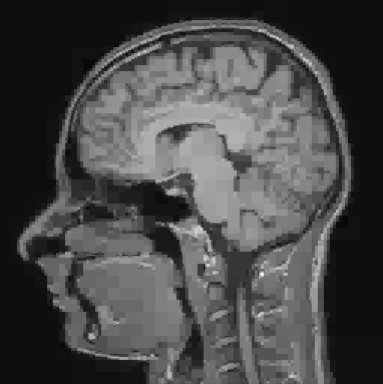}
\caption{The image on the left shows the Haar-Wavelet coefficients of a clean image, the second image shows the clean image when the only the largest 10\% of the coefficients are kept. The rest is set to zero.
On the right image less than 5\% of the coefficients are kept.}
\label{fig:wavelet_dec}
\end{center}
\vspace{-2pt}
\end{figure}
\begin{figure}[h]
\vspace{2pt}
\begin{center}
\includegraphics[width=0.27\textwidth]{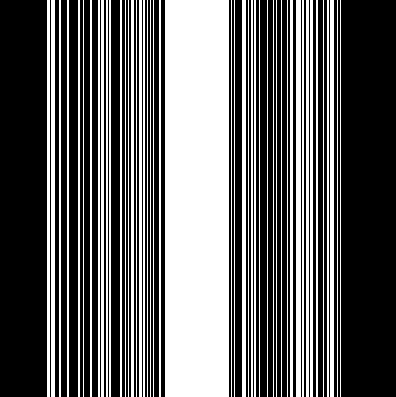}
\hspace{3pt}
\includegraphics[width=0.27\textwidth]{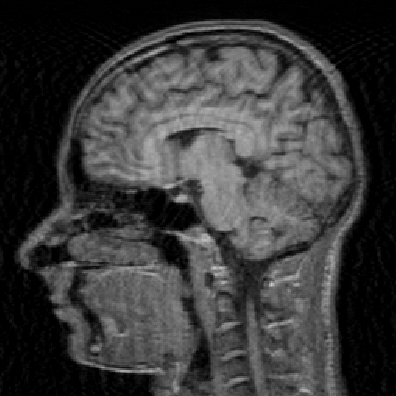}
\hspace{3pt}
\includegraphics[width=0.27\textwidth]{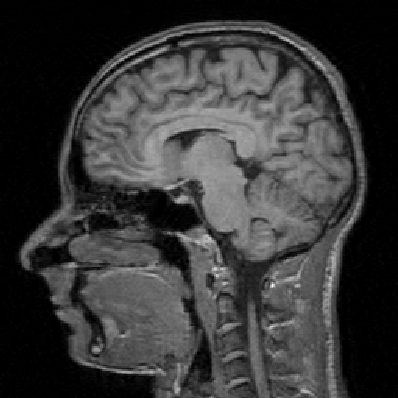}
\caption{From left to right: Sampling mask, zero filling solution, dictionary learning solution as described in Algorithm \ref{alg:ortho_dict_learn} .}
\label{fig:qual_dict_learning_results}
\end{center}
\vspace{-2pt}
\end{figure}
\begin{figure}[b]
	\centering
\includegraphics[width=0.95\textwidth]{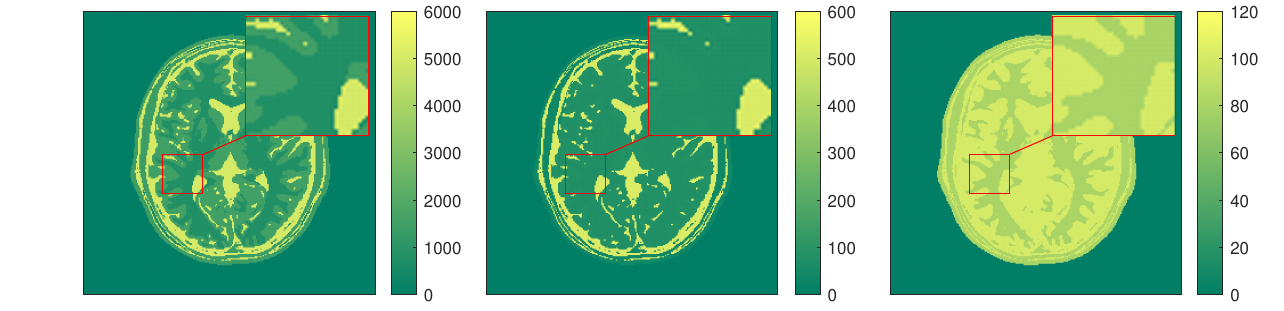} \\ \vspace{1em}
\includegraphics[width=0.95\textwidth]{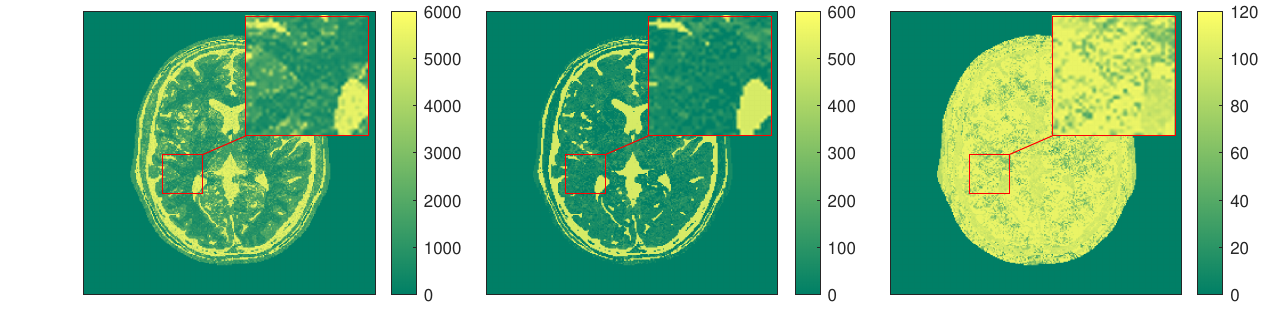} \\ \vspace{1em}
\includegraphics[width=0.95\textwidth]{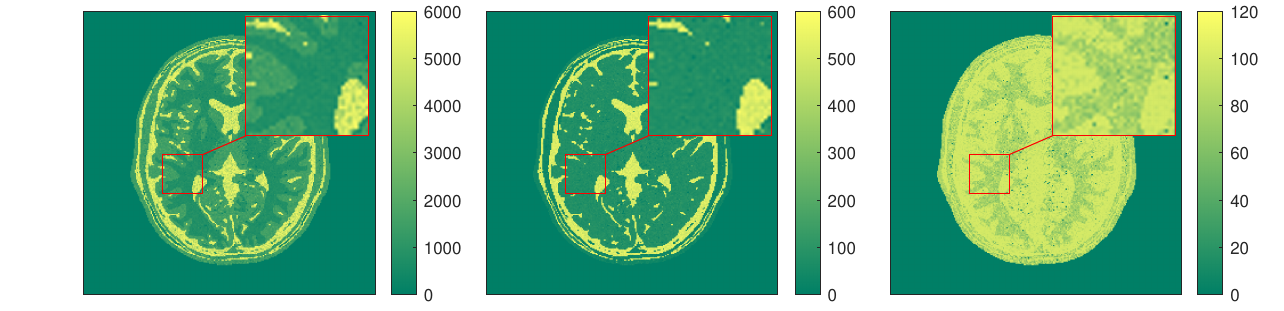} \\ \vspace{1em}
\includegraphics[width=0.95\textwidth]{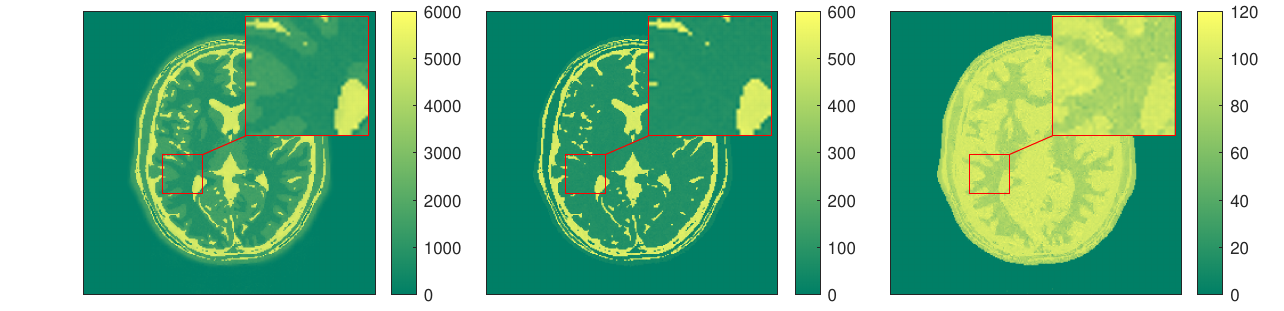} 
	\caption{From left to right: $T_{1}$, $T_2$, $\rho$ reconstruction. From top to bottom: Ground truth, BLIP, Levenberg-Marquardt, Dictionary learning approach.}
	\label{fig:qMRI_dictionary}
\end{figure}

Instead of learning the, possibly spatially dependent regularization parameter as proposed in Section \ref{sec:variational_methods}, the goal of dictionary-learning is to learn a sparse representation of a class of clean images. This is motivated by the success of Wavelet based methods and the observation that most natural images admit a sparse representation in some orthonormal basis, cf. \cite{mallat1999wavelet} for an overview on wavelets. In Figure \ref{fig:wavelet_dec} a clean image is shown together with its wavelet-decomposition. In the  second and third image a large number of wavelet coefficients  are gradually set to zero.  Still, the image retains most of its structural features which indicates that the essential information is stored in only a small number of coefficients. Since learning a basis or dictionary for the entire image has too many degrees of freedom, the classical approach is  to learn sparse representations for small image patches  \cite{aharon2006k}. During this subsection, we will work with pre-discretized complex images $u \in \mathbb{C}^{n \times m} =: X$, where $m \times n$ is the number of pixels. For such an image let $R_{i}: \mathbb{C}^{n \times m} \to \mathbb{C}^{P}$ denote the linear operator that cuts out a quadratic image patch of size $p \times p$ and puts it into a row-vector of size $P=p^2$. Here the index set $I$ is an enumeration of a selection of the possible $n \times m$ patches that can be cut out of the image $u$, see \cite{ravishankar2015efficient} for details. For simplicity we will work with $I$ indexing the set of all possible $m \cdot n$ patches of an image.  The underlying assumption is that every patch $R_{i} u_{true}, i \in I$, which is cut out of the ground-truth image can be represented as a sparse linear combination of $K \in \mathbb{N}$ image atoms $\varphi_1, \ldots, \varphi_K \in \mathbb{C}^{P}$ in the sense that
\begin{equation}
R_{i} u_{true} = \sum_{l = 1}^K c_{li} \varphi_l 
\quad \quad \|c_{i} \|_0 \leq s_{i} \ll \min\{K,P\} , \; i \in I, \label{eq:assumption_dict}
\end{equation}
where the sparsity $s_{i} \in \mathbb{N}$ of every patch $R_{i}u_{true}$ is a priori unknown but significantly smaller then the number $P$ of pixels in each patch and the number $K$ of dictionary elements. Note that \eqref{eq:assumption_dict} can also be written in a matrix vector form $Ru = DC$ where the operator $R:\mathbb{C}^{m \times n} \to \mathbb{C}^{P \times mn}$ applies $R_{i}$ to every patch indexed by $i \in I$  , $D = ( \varphi_1, \ldots, \varphi_K ) \in \mathbb{C}^{P \times K}$ denotes the unknown dictionary and $C = (c_1 , \ldots,c_{mn}) \in \mathbb{C}^{K \times mn}$ is the collection of sparse coefficients. If the number of atoms $K$ is strictly  larger then $P$ the dictionary is called overcomplete. The approach in \cite{aharon2006k} suggests to start with a training set $X = (x_1, \ldots x_M) \in \mathbb{C}^{P \times M}$ of $M$ training patches of clean images that represent a large variety of image structures. These patches are collected in a large matrix $X \in \mathbb{C}^{P \times M}$ for which the following problem is solved
\begin{equation}
\min_{D \in \mathcal{D}, C \in \mathbb{C}^{M \times K} } \frac{1}{2} \| X - DC \|^2_F  + \lambda \| C \|_s . \label{eq:KSVD_problem1}
\end{equation} 
Here $\mathcal{D} \subset \mathbb{R}^{P \times K}$ is the set of admissible dictionaries, often chosen to be the set of column normalized dictionaries, i.e. $\| \varphi_i \|_2 = 1 $ for $i = 1,\ldots, K$ and $s \in [0,1]$ defines the sparsity promoting regularizer. Note that the problem \eqref{eq:KSVD_problem1} is non-convex and non-differentiable and finding a global solution of this problem is known to be NP-hard, cf. \cite{tillmann2014computational}. Therefore possible solution algorithms focus on computing stationary points under usage of techniques from variational analysis \cite{ravishankar2015efficient, Sirotenko_PhD} or on  greedy type algorithms, see \cite{aharon2006k}. Once the dictionary $D \in \mathcal{D}$ is found, the following problem is solved
\begin{equation}
 \min_{u \in X, C \in \mathbb{C}^{P \times mn}} \frac{\mu}{2 } \| A u - y \|^2_2  + \frac{1}{2}  \| R u - D C \|_F^2 + \lambda \| C \|_s \label{eq:dict_1}
\end{equation} 
to find the estimation of $u_{true}$. For $s = 1$ the problem is convex and can be efficiently solved using classical (accelerated) splitting algorithms like FISTA \cite{beck2017first}. For $s < 1$ again the problem is non-convex (even locally non-Lipschitz for $C = 0$ which challenges solution algorithms and stationarity characterizations \cite{hintermueller2013nonconvex}). \paragraph{Blind compressed sensing}
A related method is \emph{blind compressed sensing} or \emph{ blind transform learning}, cf. \cite{ravishankar2010mr, ravishankar2015efficient, Sirotenko_PhD}. While the dictionary learning approach presented above requires the solution of the non-convex problem \eqref{eq:KSVD_problem1} in a first offline-training-phase requiring lots of training data, the authors in \cite{gleichman2011blind}  noticed that very good results can already be obtained by assuming a relatively small (often quadratic) dictionary size and by learning the dictionary simultaneously to the reconstruction process. This idea leads to the optimization problem  
\begin{equation}
 \min_{u \in X,D \in \mathcal{D}, C \in \mathbb{C}^{P \times mn}} \frac{\mu}{2} \| A u - y \|^2_2  + \| R u - D C \|_F^2 + \lambda \|C \|_s. \label{eq:general_linear_dict_problem}
\end{equation} 
Let us briefly sketch the details of  \cite{ravishankar2015efficient} where an efficient update strategy is combined with a convergence analysis. We will focus on learning an orthogonal transform here, while the approach presented in \cite{ravishankar2015efficient} also allows for more general dictionaries. This corresponds to the choice 
\begin{equation}
\mathcal{D} = \{ D \in \mathbb{C}^{K \times K} \, | \, D^H D = I\}. 
\end{equation}
The general idea of the solution algorithm is to update one variable while keeping the others fixed and then iterate cyclically.  We start with the update of the coefficients, i.e. we want to optimize over $C \in \mathbb{C}^{K \times K}$ while keeping $D$ and $u$ fixed. The resulting optimization problem in the $k$-th step of the algorithm reads for given $D_k \in \mathcal{D}$ and $u_k \in X$:
\begin{align}
C_{k+1} \in   \argmin_{C \in \mathbb{C}^{P \times mn}}\; \frac{1}{2}\| R u_k - D_k C_k \|_F^2 + \frac{\lambda_{C,k}}{2}\|C - C_k\|_F^2 + \lambda \|C\|_s ,  \label{eq:C_update}
\end{align}
whose solution has a closed form in terms of the well known proximal operator, see \cite{rockafellar2009variational}. Here we took advantage of the orthogonality of the transform $D^{-1} = D^H$.  Let us now consider the update for the transform $D$ for fixed $u_k$, $C_{k+1}$ and previous iterate $D_k$. We want to calculate in this step
\begin{equation}
D_{k+1} \in  \argmin_{D \in \mathcal{D}}\; \frac{1}{2} \|D C_{k+1} - Ru_k \|^2_F + \frac{\lambda_{D,k}}{2}\|D - D_k \|^2_F. \label{eq:dictionary_update}
\end{equation}
It is remarkable that the non-convex problem in \eqref{eq:dictionary_update} which involves the solution of a quadratic objective on the matrix manifold $\mathcal{D} \subset \mathbb{C}^{K \times K}$ has a closed form solution, that can be computed as fast as the singular-value-decomposition of $Ru_k C^H + \lambda_{D,k} D_k$ can be computed. In general the minimizer is not unique, cf. \cite{ravishankar2015efficient}. Once we have updated the dictionary $D$ the last step requires the update for the image $u$. Therefore the following problem needs to be solved
\begin{equation}
u_{k+1} \in \argmin_{u}\; \frac{\mu}{2} \|A u - y \|_2^2 + \frac{1}{2} \| Ru - D_{k+1}C_{k+1} \|^2_F + \frac{\lambda_{u,k}}{2} \|u - u_k \|^2. 
\end{equation}
The equivalent first-order system of this quadratic optimization problem reads
\begin{equation}
\left( \mu A^H A + R^\top R + \lambda_{u,k} \mathrm{Id} \right) u = A^H y + R^\top D_{k+1}C_{k+1} + \lambda_{u,k} D_{k+1} \label{eq:general_u_update}
\end{equation}
The system in \eqref{eq:general_u_update} is in general very large and can only be solved approximately for a general forward operator $A$, e.g. using methods of conjugate gradient type  \cite{golub2013matrix}. However in the case of the MRI operator and a special form of $R$ this system has a closed form solution which can be computed easily \cite{ravishankar2010mr,
ravishankar2015efficient} for details. The schematic overall algorithm for the linear MRI problem is depicted in Algorithm \ref{alg:ortho_dict_learn}.
\begin{algorithm}[h]
\begin{algorithmic}
\State Given $u_0 \in X, D_0 \in \mathcal{D}, C_0 \in \mathbb{C}^{P \times nm}$ and a sequence $\{\lambda^u_n,\lambda^D_n,\lambda^C_n\}_{n\in\mathbb{N}}$ of positive real numbers, iteratively solve the following problems for $n=0,1,2,\ldots$:\\
\begin{itemize}
\item[(1)] Compute $C_{k+1}$ given $u_k,D_k$  by formula \eqref{eq:C_update}.    
\item[(2)] Compute $D_{k+1}$ formula \eqref{eq:dictionary_update}.
\item[(3)] Compute $u_{k+1}$ using formula \eqref{eq:general_u_update}.
\end{itemize}
\State If a stopping criterion is met, return $u_k$.  
	\caption{\textbf{Blind Compressed Sensing \cite{ravishankar2010mr}}}\label{alg:ortho_dict_learn}
\end{algorithmic}
\end{algorithm}
The procedure can be shown to converge against a stationary point in the limiting sense \cite{rockafellar2009variational} under the assumption that the overall objective function satisfies the so called Kurdyka-Lojasiewicz-inequality (KL-inequality)  \cite{attouch2013convergence, attouch2010proximal} and there exists a constant $c>0$ such that $\min \{ \lambda_{k,u},\lambda_{k,D},\lambda_{k,C} \} \geq c $ for every $k $ . In \cite{ravishankar2015efficient} a refined convergence result is presented which does neither require the KL-inequality nor the step-size bounds but only yields convergence up to a subsequence. For a proof of this result and details see again \cite{ravishankar2015efficient}. A prototypical numerical result of the discussed approach is presented in Figure \ref{fig:qual_dict_learning_results}.
\paragraph{Blind compressed sensing for qMRI} 
Despite the relatively high computational effort the methodology of blind dictionary learning has recently been proven to be successful also for qMRI, \cite{kofler2023quantitative, Sirotenko_PhD}. We recall the optimization problem for qMRI which was considered also in \eqref{eq:Re_C_Opt}:	
\begin{equation*}
	\underset{q\in  \mathcal{C}_{ad}}{\operatorname{min} }\; \frac{\mu}{2} \| A\circ \Pi(q)-y \|_{2}^2 + \frac{1}{2} \| q \|^2_{2}, 
\end{equation*}
with the minor modification that we put the penalty parameter in front of the data discrepancy term. However instead of using a Tikhonov  regularisation approach with a smooth quadratic penalty as proposed in \cite{DonHinPap19}, a regularization strategy via blind dictionary learning is proposed in \cite{kofler2023quantitative} and some ongoing work \cite{Sirotenko_PhD}. Analogously to the corresponding method for qualitative MRI which is described at the beginning of this subsection we aim at solving the optimization problem
\begin{equation}
	\underset{q \in  \mathcal{C}_{ad}}{\operatorname{min} }\; \frac{\mu}{2} \| A\circ \Pi(q)-y \|_{2}^2 + \mathcal{R}(q) ,
\end{equation}
where the regulariser $\mathcal{R}(q)$ is given by another minimization problem
\begin{equation}
\mathcal{R}(q) = \inf_{D \in \mathcal{D},C \in \mathcal{C}} \frac{1}{2} \|Rq - DC \|_F^2 + \alpha \|C\|_s.
\end{equation} 
Here again $s\in [0,1]$ specifies the sparsity promoting regularizer which was also proposed in the linear case and the set $\mathcal{C}$ is the set of sparse coefficients. $R$ denotes the linear patch-extraction operator and $\mu, \alpha >0$ balance overall the regularization effect introduced by learning the dictionary. The eventual optimization problem reads
\begin{equation}
	\underset{(q,D,C) \in  \mathcal{C}_{ad} \times \mathcal{D} \times \mathcal{C}}{\operatorname{min} }\; \frac{\mu}{2} \| A\circ \Pi(q)-y \|_{2}^2 +  \frac{1}{2} \|Rq - DC \|_F^2 + \alpha \|C\|_s \label{eq:dict_qmri}
\end{equation}
In \cite{Sirotenko_PhD} a model based optimization algorithm is proposed in order to find limiting-stationary points of \eqref{eq:dict_qmri}. While the update steps for the sparse coefficients $C$ and the dictionary $D$ look similar to the steps in \eqref{eq:C_update} and \eqref{eq:dictionary_update} respectively, the update step of the tissue parameter $q$ does not require the solution of a linear equation as in \eqref{eq:general_u_update} but the solution of a non convex constrained optimization problem. For this purpose a Levenberg-Marquardt approximation is employed and the eventual update step reads
\begin{equation}
q_{k+1} = \argmin_{q \in \mathcal{C}_{ad}} \;\frac{\mu}{2} \| A \circ \Pi' (q_k)[q - q_k]+  A \circ \Pi (q_k) - y \|_2^2 + \| Rq - D_{k+1}C_{k+1} \|_F^2 + \frac{\lambda_{q,k}}{2} \|q - q_k\|_{2}^2. 
\end{equation} 
In \cite{Sirotenko_PhD} a line-search strategy is proposed to find the step size-parameter $\lambda_{q,k}>0$ that has to be adjusted in dependence of the unknown Lipschitz constant of $\Pi'$ in order to guarantee a sufficient descent of the objective in each step of the algorithm. The overall method alternates between the three optimization steps and is shown to converge against a limiting-stationary point. Prototypical numerical results from \cite{Sirotenko_PhD} are depicted in Figure \ref{fig:qMRI_dictionary}, where the dictionary learning approach is compared against the pure Levenberg-Marquardt method from Section \eqref{sec:quantitative_MRI} and the BLIP method in \cite{DavPuyVanWia14}. Here the number of image frames was chosen to be $L=200$ and the under-sampling factor is 16. Moreover small complex noise was added to the under sampled Fourier data. The details can be found in \cite{Sirotenko_PhD} where also extensions of the methods are discussed. 

\begin{figure}[b]
	\centering
\includegraphics[width=0.95\textwidth]{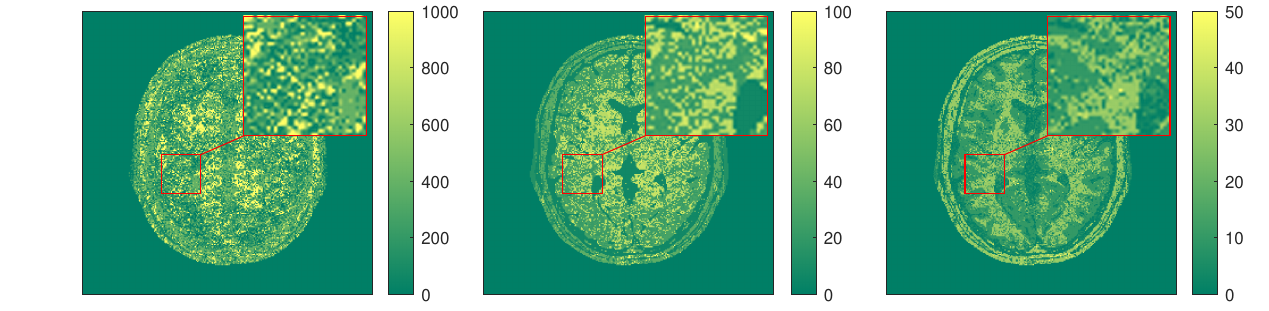} \\ \vspace{1em}
\includegraphics[width=0.95\textwidth]{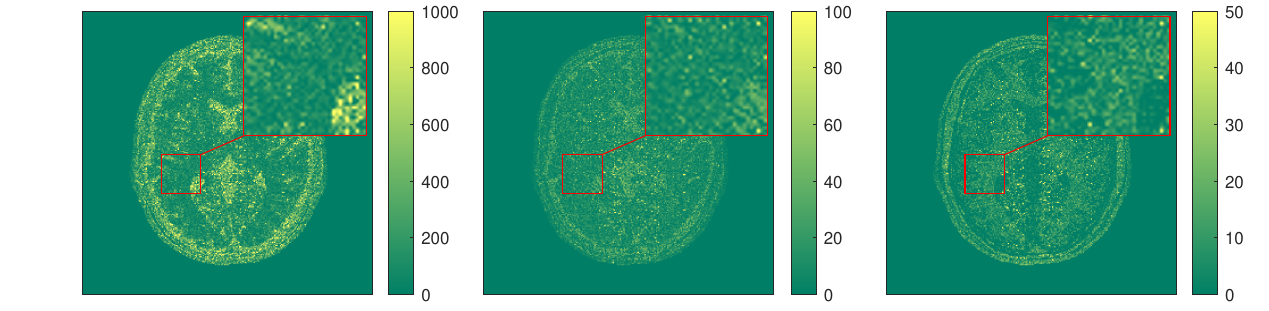} \\ \vspace{1em}
\includegraphics[width=0.95\textwidth]{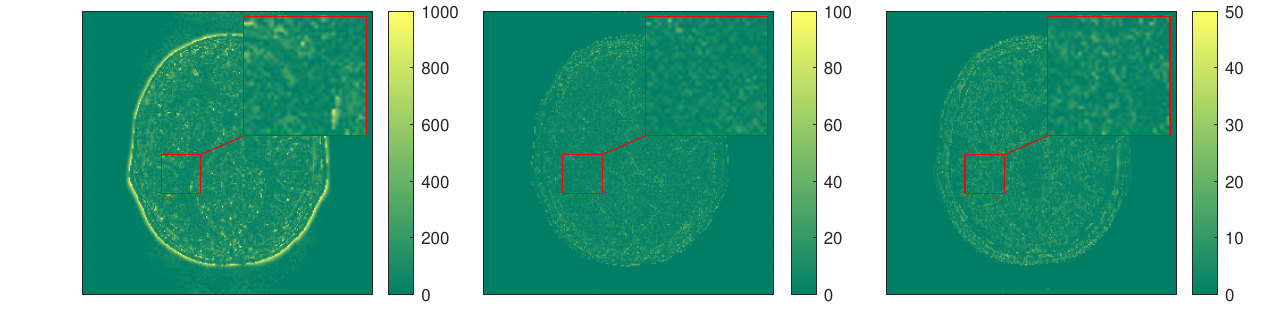} 
	\caption{From left to right: Error-maps for the $T_1$, $T_2$ and $\rho$ reconstruction. Top to bottom: BLIP, Levenberg-Marquardt, Dictionary learning approach.}
	\label{fig:qMRI_dictionary_errors}
\end{figure}
Despite the practical success of the dictionary learning methods described above, there are still many open theoretical questions. First of all even for qualitative MRI the corresponding optimization problem  \eqref{eq:general_linear_dict_problem} is a non-convex and non-differentiable program  which depends on the initialization, the optimization algorithm, the parameter choice, and also the choice of the sampling pattern. Moreover it is a priori unclear if a solution of \eqref{eq:general_linear_dict_problem} or even more difficult a solution of \eqref{eq:dict_qmri} has anything to do with ground truth  data-generating dictionary. This question of the so called identifiability has been discussed in \cite{spielman2012exact,geng2014local,schnass2015local,gribonval2015sparse}    
but is in general still an open problem. Further theoretical progress has concentrated mostly on different facets of the classical dictionary learning problem in \eqref{eq:KSVD_problem1}. In \cite{sun2016complete, sun2016complete2} the structure of minimizers and stationary points is studied. More sophisticated minimization algorithms are for instance considered in e.g.  \cite{khanh2022block, li2022riemannian}. The question of whether alternating minimization schemes are able to recover the data generating dictionary is investigated in \cite{arora2015simple, agarwal2016learning, ravishankar2020analysis, liang2022simple} under some specific probabilistic model. 
\begin{tcolorbox}[
enhanced jigsaw,drop shadow, colback= yellow!75!black, boxsep=0.1cm,
boxrule=1pt, width=1\textwidth, 
interior style={top color=mygray!20!white,
bottom color=mygray!20!white}, 
 opacityback=1,
fonttitle=\bfseries, arc=5pt]
\begin{remark}[Function space versions of dictionary learning] 
Most of the work on patch-based dictionary learning was done in the finite dimensional set-up, i.e. a pre-discretized image with a fixed number of pixels was assumed. 
Generalisations of patch-based  dictionary learning are inherently difficult as the generalisation of the patch extraction operator is not clear. The straightforward way would  first apply a discretization to an image in the function space $X$ and cut the patches out of the discretization. But in this way $R$ has a nullspace where no regularization acts at all and also the solution again depends on a pre selected discretization. Recently a variant of dictionary learning called \emph{convolutional dictionary learning} was proposed which does not rely on patches but on kernels which are convolved against sparse coefficients
\cite{garcia2018convolutional,
papyan2017working} to reconstruct the image. Convolutional dictionary learning admits a natural extension to infinite dimensions which is studied in \cite{chambolle2020convex}. In addition to these difficulties on the modeling side, also the extension of the convergence theory, in particular of the KL-inequality,  towards infinite dimensions seems delicate. While there are some generalizations towards Hilbert spaces \cite{bolte2010characterizations} and also to general metric spaces \cite{hauer2019kurdyka} available which are used to study the convergence of gradient flow problems, the application to inverse problems is rarely studied in the literature. This is due to the fact that the verification of the KL-inequality is more involved since its connection  to real-semialgebraic geometry, as it is usually used for finite dimensional problems \cite{attouch2013convergence, attouch2010proximal}, is to the best of our knowledge not available.
\end{remark}
\end{tcolorbox}

\subsection{Learning-informed physics in quantitative imaging}\label{sec:learning_physics}

\subsubsection{General set-up}
The quantitative methods for imaging discussed in Section \ref{sec:quantitative_MRI} rely on the assumption that the equations that govern the physical constraints of the inverse problem are known a priori.
This is a valid assumption if the governing physical model is derived from well-established physical arguments that are validated experimentally.
Often, however, the physical equations are of phenomenological nature only and rooted in empiricism,
such that it might be advisable to substitute the physical model or parts of it by learning-informed components, which are trained based on available physical data.
\emph{Learning-informed models} have the advantage that they are versatile in describing correlations and structures in given physical data.
Consequently, utilizing learning-informed models reduces the risk of introducing modeling errors caused by false physical assumptions.
Further, learning-informed models can reduce the computational cost associated with the evaluation of the models, e.g., when learning an efficient-to-evaluate control-to-state map for control problems like those described in Section \ref{sec:quantitative_MRI}.
In the following, we discuss some recent advances of learning-informed physical constraints for general inverse problems and then in particular with respect to the qMRI problem.

\begin{figure}
	\centerline{
		\includegraphics[width=0.90\textwidth]{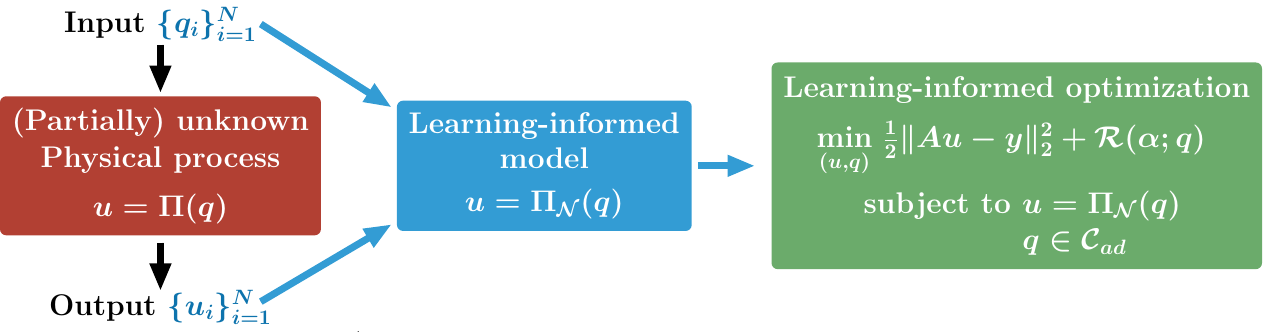} \hspace{0.3cm}
	}
	\caption{Learning a (partially) unknown physical process $\Pi$ and embedding the learned map $\Pi_{\mathcal{N}}$ to the integrated-physics optimization problem.}
	\label{fig:diagram}
\end{figure}

We begin by considering a  physics-constrained inverse problem, which is a general version of the integrated-physics qMRI approach \eqref{eq:C_Opt} discussed in Section \ref{sec:quantitative_MRI}
\begin{equation}
\left \{
\begin{aligned}
\label{eqmf:constrained_inverse_problem}
&\min_{u,q\in \mathcal{C}_{ad}} \;\frac{1}{2} \| Au - y \|_{2}^2 + \mathcal{R}(\alpha ; q),\\
&\text{subject to }\;\;  e(u,q) = 0.
\end{aligned}\right.
\end{equation}
Recall here that  the (possibly nonlinear) physical constraint $e(u,q) = 0$, which may take the form of an ordinary or partial differential equation, relates the state variable $u$ and the control variable $q$.
We assume that the physical constraint $e(u,q) = 0$ is well-posed such that for a given control $q$, the state $u$ that satisfies the physical constraint is unique, i.e., we can define (e.g., via the implicit function theorem) the control-to-state (or parameter-to-solution) map
\begin{equation}
u = \Pi(q).
\end{equation}
In  Section \ref{sec:quantitative_MRI}, it was assumed that the physical constraint $e(u,q) = 0$ and the control-to-state operator $\Pi$ are a priori known, and in particular it could be expressed by the Bloch equations.
As aforementioned, in this section, we seek to substitute the physical model or parts of it by learning-informed components.
Here we distinguish between two conceptually different strategies.
The first strategy is to learn the physical model in its implicit form $e(u,q) = 0$ from data.
If $e(u,q) = 0$ is an ordinary or partial differential equation, this could mean to learn the unknown differential operator or an unknown nonlinear term in the differential equation.
This strategy has been followed by \cite{schmidt_distilling_2009,brunton_discovering_2016,rudy_data-driven_2017,raissi_deep_2018,raissi_hidden_2018,chen_physics-informed_2021, DonHinPap22, DonHinPap22b, DonHinPapVoe22}, to name only a few.
When the target of learning is the implicit form of the physical model, computing the state for a given control is still needed after the learning has taken place and it usually requires solving a differential equation.
The alternative, conceptually different strategy is to directly learn the explicit control-to-state operator $\Pi$ from available data such that the computation of the state for a given control becomes computationally more efficient and closer to the true physical process.  This approach was introduced  in \cite{DonHinPap22} and the overall idea is depicted in the diagram in Figure \ref{fig:diagram},
The focus of the following discussions lies on this latter strategy because it is well-aligned with the integrated-physics approached discussed in Section \ref{sec:quantitative_MRI}.
The  strategy for substituting the physical model by a learning-informed ansatz is to parameterize the unknown control-to-state operator, i.e.,
\begin{equation}
\label{eq:control_to_state_parameterization}
u = \Pi(q) \approx \Pi_{\mathcal{N}}(\theta;q),
\end{equation}
where $\theta\in\mathbb{R}^\ell$ is a set of $\ell$ tunable parameters, for example the weights of a neural network. If labeled data pairs $(q_i,u_i)$ with $i\in\{1,\dots,N\}$ are available, the tunable parameters can be calibrated by minimizing the mismatch between model prediction and data in a supervised learning fashion
\begin{equation}
\label{eq:training}
\min_{\theta} \frac{1}{N}\sum_{i=1}^N\| \Pi(\theta;q_i) - u_i \|_{2}^2 + r(\theta).
\end{equation}
where $r$ denotes some potential regularization for $\theta$ reflecting prior information which one wants to impose on $\theta$. Also, in some cases additional constraints on $\theta$ might be relevant. Solving \eqref{eq:training} is often a very challenging tasks, giving rise to many current efforts in solver design. Also, non-smooth variants of \eqref{eq:training} are relevant in practice.

\subsubsection{Learning the Bloch solution map for qMRI}
\label{sec:fin_dim_learning}

We turn back to the physics-constrained inverse problem in \eqref{eqmf:constrained_inverse_problem} and consider it in the context of qMRI, i.e., we consider the Bloch equations as the physical constraint (see Section \ref{sec:quantitative_MRI}) where, following the terminology from optimal control, the state variable is the magnetization and the control variables are the physical tissue parameters.
The literature reports different attempts to tackle the Bloch equations with machine learning. 
These attempts can be classified into approaches that seek to learn the Bloch solution map or its inverse, or both through encoder-decoder type models, as discussed in the following.

Learning the inverse of the Bloch solution map can, for example, be used to efficiently identify the physical tissue parameters given the magnetization response,
by leveraging deep neural networks as demonstrated by \cite{hoppe_deep_2017,virtue_better_2017,cohen_mr_2018,oksuz_magnetic_2018,golbabaee_geometry_2019,hoppe_magnetic_2019}.
For instance, the work in \cite{cohen_mr_2018}  introduced the framework DRONE (Deep RecOnstruction NEtwork), in which the inverse Bloch solution map is approximated by a neural network that is trained on a simulated MRF look-up-table.
The trained network requires less memory and allows for faster predictions than the MRF look-up-table.
In the other direction, neural networks can for example be used to efficiently generate the Bloch look-up-table for MRF as shown by \cite{yang_game_2020} using generative adversarial neural networks (GANs), by \cite{hamilton_machine_2020} using feed-forward ANNs, and by \cite{liu_fast_2021} using recurrent neural networks (RNNs).

Another approach in this context, aligned with Figure \ref{fig:diagram},  is to learn the Bloch solution map and substitute it in the  reduced (integrated physics) formulation given in \eqref{eq:Re_C_Opt}.
This idea was realized in \cite{DonHinPap22} where feed-forward ANNs have been utilized for learning the control-to-state operator $\mathcal{N}(q) \approx \Pi(q) $. This resulted in the following learning-informed optimization problem:
\begin{equation}\label{eq:NN_model}
\underset{q\in  \mathcal{C}_{ad}}{\operatorname{min} }\; \frac{1}{2} \| A\circ \mathcal{N}(q)-y\|_{2}^2 + \frac{\alpha}{2} \| \nabla q \|_{2}^2, 
\end{equation}
where for simplicity a standard Tikhonov regularisation was employed for $q$ albeit with small $\alpha$ in order to stabilize the solution process.
In Algorithm \ref{alg:learning}, we summarize the above procedure. We note that in this approach the neural network $\mathcal{N}$ is trained in an \emph{offline phase}, and then embedded into the minimization problem \eqref{eq:NN_model}. However, integrated learning techniques, i.e., training $\mathcal{N}$ while solving \eqref{eq:NN_model}, are conceivable.

\begin{algorithm}
\begin{itemize}
\item[(1)] Collect a set of training data which are resulting from the targeting physical models with (experimental or numerical) measurements.
\item[(2)] Build a neural network architecture to learn this physical model using the training data from Step (1). 
\item[(3)] Solve the optimization problem \eqref{eq:NN_model} with learning informed operators or constraints with learning informed differential equations.
\end{itemize}
	\caption{\text{A sketch of the learning-informed models for quantitative imaging}}
	\label{alg:learning}
\end{algorithm}

\begin{remark}
The description of Algorithm \ref{alg:learning} looks quite abstract. Here we take an example to explain the ideas with more details.
For our qMRI example, the training data are collected from the dictionary which is introduced in the MRF paper \cite{Ma_etal13}, see also \cite{DonHinPap19}. The data consist of time series which mimic the discretization of the Bloch equations, where the input is the tissue parameter, and the output is the magnetization. Using this training data, we can learn the map from tissue parameters to magnetization (it is a type of Nemytskii operator in the case of \cite{Ma_etal13}). However, the learning approach can be more general in the case where the Bloch equations are considered to be a simplification of a more complicated physical process, and more  data are available for the training.
The architecture of the neural network then depends on the format of the training data. For the time series data, residual type neural networks are popular, see a description in \cite{DonHinPap22}, which is also the one we have used in the numerical examples of this review.
In terms of the third step, there are typically many ways to solve the optimization problem. In the considered example, we have used a semismooth Newton type algorithm \cite{HIK} to take care of the box constraints for the tissue parameters. For more details one can also refer to \cite{DonHinPap22}.
\end{remark}

\begin{figure}[h]
	\centering
\includegraphics[width=0.95\textwidth]{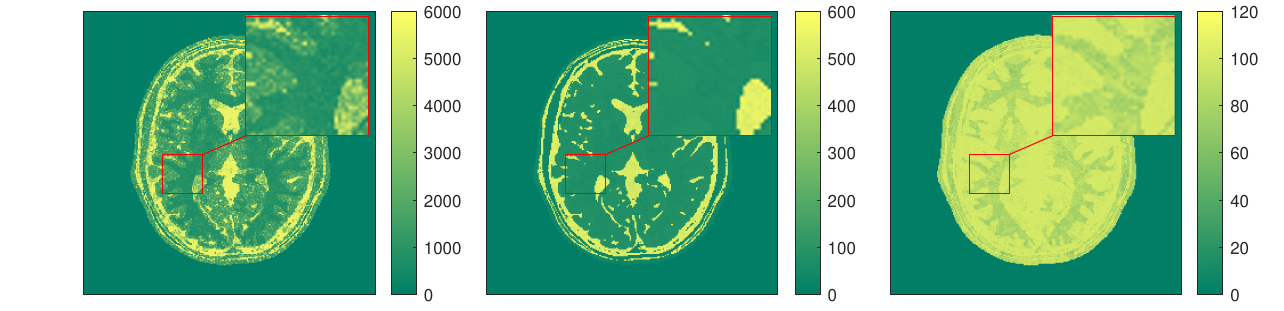}\\
\includegraphics[width=0.95\textwidth]{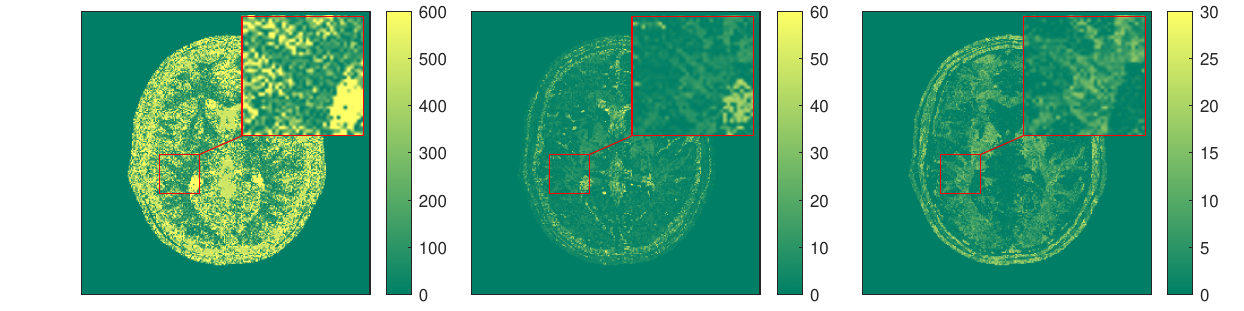} 
	\caption{The estimated parameters $T_{1}$, $T_{2}$ and $\rho$ using the optimization framework with learning-informed models.  Top row are the estimated parameters, and bottom row are the relative errors. By comparing the result with the ones in Figures \ref{fig:qMRI_comparisons_parameters} and \ref{fig:qMRI_comparisons_errors}, we see that using learning-informed operator, the obtained results are comparable, and the computational time is less in solving the optimization problem, since the learning is done before hand, and heavy evaluation of the equation is avoided in every iteration.
	}
	\label{fig:qMRI_NN}
\end{figure}

We present some numerical examples in Figure \ref{fig:qMRI_NN}, where we compute solutions of the optimization problem \eqref{eq:NN_model}. In this example, the neural network simulating the Bloch map is trained via a training set which was obtained through discrete dynamics of simulated Bloch data. In particular, a set of parameter pairs $\{q_{i}\}_{i=1}^{N}$ in the feasible domain of $T_1$ and $T_2$ is taken as inputs of the network, and then the simulated Bloch dynamics resulting from all these pairs are regarded as the  outputs. The simulated k-space data in this example are generated using exactly the same setting as the one used for  Figures \ref{fig:qMRI_comparisons_parameters} and  \ref{fig:qMRI_comparisons_errors}. 
 For such a discrete time series, the Bloch map is a Nemytskii type operator, and thus can be approximated sufficiently well in a uniformly bounded feasible set. 

\begin{tcolorbox}[
enhanced jigsaw,drop shadow, colback= yellow!75!black, boxsep=0.1cm,
boxrule=1pt, width=1\textwidth, 
interior style={top color=mygray!20!white,
bottom color=mygray!20!white}, 
 opacityback=1,
fonttitle=\bfseries, arc=5pt]
\begin{remark}[Learning-informed qMRI in function space]
The learning-informed optimization problem for qMRI \eqref{eq:NN_model} can be also posed and analyzed in a function space setting \cite{DonHinPap22}. The key observation here is that the Bloch solution map $\Pi$ and its approximation $\mathcal{N}$ can be regarded as Nemytskii operators, with a voxel-wise action, i.e., $\hat{\Pi}, \hat{\mathcal{N}}:\mathcal{C}_{ad} \to [(L^{\infty}(\Omega))^{3}]^{L}$: 
\begin{align}
\hat{\Pi}(\hat{\rho}, \hat{T}_{1}, \hat{T}_{2})(x)&= \Pi(\rho(x), T_{1}(x), T_{2}(x))\label{Pi_nem}\\
\hat{\mathcal{N}}(\hat{\rho}, \hat{T}_{1}, \hat{T}_{2})(x)&= \mathcal{N}(\rho(x), T_{1}(x), T_{2}(x))\label{N_nem}
\end{align}
for almost every $x\in\Omega$.
Here $\mathcal{C}_{ad}$ is a subset of $[L^{\infty}(\Omega)^+]^{3}$, i.e. $L^{\infty}$-functions that are bounded away from zero, denoting feasible values for $\hat{\rho}, \hat{T}_{1}, \hat{T}_{2}\in L^{\infty}(\Omega)$. Note that $\hat{\mathcal{N}}$ can be in principle learned via the network $\mathcal{N} \approx \Pi$ with inputs in $\mathbb{R}^{3}$. Under this setting, the analysis of \eqref{eq:NN_model} in the function space setting can be performed. For instance in \cite{DonHinPap22} convergence of the solution of the learning-informed problem  to the solution of the problem given in \eqref{eq:Re_C_Opt} was proved upon increasing the approximation capability of the low dimensional network $\mathcal{N}$ and corresponding error estimates were also provided.\\[0.5em]
\emph{Towards learning infinite dimensional physics operators for qMRI}. Conventionally, the  universal approximation theorem \cite{hornik_multilayer_1989} considers neural network functions between finite dimensional spaces.
The idea of studying neural networks as approximators of functionals on infinite dimensional spaces 
has been initiated by \cite{chen_universal_1995}, where two  architectures were introduced; one for approximating functionals (see also the related earlier work \cite{sandberg_approximations_1992}) and one for approximating operators.
The number of parameters in these models
 can be substantially larger than in a standard feed-forward network.
Thus, training such a model requires more computational resources; and \cite{chen_universal_1995} did not provide a practical training example.
In recent years, as computers and optimizers for machine learning models have become more and more powerful, a blossoming of the idea of learning functionals or operators is observed \cite{lu_learning_2021,bhattacharya_model_2021,li_fourier_2021,kovachki_universal_2021,raonic_convolutional_2023}.
Prominent ANN architectures for learning operators between infinite-dimensional spaces are DeepONets \cite{lu_learning_2021}, which are inspired by the architecture presented in  \cite{chen_universal_1995}, and (Fourier) Neural Operators.
We see the opportunities of operator learning within the problems setting of the current paper, where the qMRI physical model can be realized as an operator as well. Once a mature operator learning scheme has been developed, it will not only improve the modeling process but also contribute efficient solvers.
Note that more realistic qMRI dynamics may involving partial differential equations. In this context, operator learning is an even more recent topic and one may refer to the papers \cite{lu_learning_2021,bhattacharya_model_2021,li_fourier_2021,kovachki_universal_2021,raonic_convolutional_2023} and the references therein.
\end{remark}
\end{tcolorbox}

\section{Towards FAIR data}

Similar to all mathematical research, the imaging field is intrinsically tied to research data.
This encompasses the image data $y$ used for applying the mathematical methods, its parameter sets, and even the algorithms themselves along with intricate workflows employed for image processing.
The digitalization trend has spurred substantial data growth. 
Particularly evident in data-driven methodologies, this data's availability has become fundamental for research. 
However, overlooking effective research data management could result in ``dark data''~\cite{Schembera2020}, posing challenges to research reproducibility and contributing to the broader issue known as the reproducibility crisis in science~\cite{Baker2016}.

In response to these challenges and the overarching call for open science and open data, Tim Berners-Lee's 5-star principles~\cite{TBL5stars} and the more recent FAIR (Findable, Accessible, Interoperable, Reusable) principles~\cite{wilkinson_fair_2016} have emerged. 
These principles provide guidance and categorization of compliance levels for research output and data. 
Similarly, the complementary FAIR4RS~\cite{FAIR4RS} principles have been formulated to address the specific requirements of research software.

Despite their distinct characteristics, open data and FAIR principles share a unifying objective: enhancing the accessibility of research data, much like how journal publications have been provided to the scientific community, facilitating the practice of building upon existing knowledge. 
Open data and open-source software should be accessible to all under fitting open licenses. 
In contrast, the FAIR principles delineate the methods and conditions for human and machine data access. 
Moreover, they introduce stringent prerequisites for metadata descriptions and formats, enhancing findability, interoperability, and potential reusability.

Specifically, findability encompasses assigning a distinct identifier, furnishing comprehensive metadata descriptions, and indexing the data within a resource that's both searchable by humans and machines, streamlining researcher access. 
Accessibility ensures that researchers with interest can approach the data using standardized, open, and cost-free communication protocols. 
The metadata remains perpetually accessible, even if the data itself isn't or becomes inaccessible over time. 
Interoperability empowers data use across different contexts from their original purpose. 
Reusability necessitates meticulous documentation aligned with community standards, facilitating reproducibility and enabling other researchers to utilize a dataset for their own applications. 
This reusability also mandates releasing the data with a fitting license.

In an attempt of a community-driven approach to make research data FAIR and provide appropriate infrastructures for storage, collection of the scientific knowledge, and innovative services to enable new scientific results, the German government decided to fund the building of a National Research Data Infrastructure~\cite{nfdi}. 
In particular, the community approach is realized by discipline-specific consortia that address the needs of the respective area in science. 
The consortium for mathematics, the Mathematical Research Data Initiative (MaRDI)\cite{mardiantrag,mardi,mardiportal}, has started its work in 2021 and is gradually releasing new services to the mathematical community such as specific data repositories, services for creating and maintaining a data management plan and a MaRDI knowledge graph to connect the mathematical knowledge, software, publications, and data.

Mathematical research projects in the area of imaging may benefit from an existing data infrastructure while also contributing to its contents. 
This already starts at the proposal phase, where an appropriate research data management plan has to be developed in accordance with the requirements of the funding agency~\cite{MaRDI-DFG-WP}.
Besides the description of the data used and created within the project such a plan usually contains the elaboration on its documentation, e.g., with rich metadata, thoughts on storage in appropriate data repositories, assignment of suitable licenses, descriptions of availability to third-parties and responsibilities within the project, and more general measures to increase the FAIRness of the project data and software. 
Such a research data management plan is usually a living document adjusted to the needs of the project over time. 

Projects for mathematical imaging general should consider research data such as: images used for algorithmic evaluation, collections of images used for training and testing, software scripts for processing examples, implementations of algorithms in some more general-use software, knowledge about the algorithm and related publications and software as contribution to the MaRDI knowledge graph. 
Naturally, some of these data are generated within the project, others are often re-used from other available resources, such as OpenNeuro~\cite{openneuro}, the Human Connectome Project~\cite{HCP}, the Alzheimer’s Disease Neuroimaging Initiative~\cite{adni}, or \cite{brainweb} for simulated brain data in the field of neuroimaging.
It should be considered natural for research project to publish its generated data and software in a FAIR manner as described above in the same way as it is common for the scientific articles describing and discussing the scientific results. 
This will increase visibility of the research, enhance its reproducibility, prevent re-invention of the wheel by duplicated efforts for data, software, and knowledge creation, and lower the burden to start new research.

\subsection*{Acknowledgments}
GD acknowledges  the support from NSF of Hunan Province (project No. 2024JJ5413). \\
This research is funded by the Deutsche Forschungsgemeinschaft (DFG, German Research Foundation) under Germany's Excellence Strategy - MATH+: The Berlin Mathematics Research Center [EXC-2046/1 - project ID: 390685689].

\subsection*{Author contributions}
The various authors have made the following contributions according to the CRediT-taxonomy.\\[0.5cm]
\textbf{Guozhi Dong}: Writing-original draft, Writing–review \& editing, Software,  Conceptualization,  
\textbf{Moritz Flaschel}: Writing-original draft,
\textbf{Michael Hintermüller}: Project administration, Supervision, Writing–review \& editing, Conceptualization, Funding acquisition, 
\textbf{Kostas Papafitsoros}: Project administration, Supervision, Writing-original draft, Writing–review \& editing, Conceptualization
\textbf{Clemens Sirotenko}: Writing-original draft, Writing–review \& editing,  Software, Conceptualization, 
\textbf{Karsten Tabelow}: Writing-original draft, Writing–review \& editing, Software.

\subsection*{Financial disclosure}

None reported.

\subsection*{Conflict of interest}

The authors declare no potential conflict of interests.

\bibliographystyle{plain}
\bibliography{BibGD,BibMF,BibMH,BibKP,BibCS,BibKT}%

\end{document}